\documentclass[onefignum,onetabnum]{siamart220329}

\usepackage{lipsum}
\usepackage{amsfonts}
\usepackage{graphicx}
\usepackage{epstopdf}
\usepackage{algorithmic}
\usepackage{amsopn}

\usepackage{amsmath}
\usepackage{amssymb}
\usepackage{fancyhdr}
\usepackage{color,verbatim}

\usepackage[caption=false]{subfig}
\usepackage{multirow}
\usepackage{rotating}
\usepackage{tabularx}
\usepackage{ltablex}
\usepackage{booktabs}
\usepackage{listings}
 \allowdisplaybreaks[4]
\usepackage{longtable}

\newsiamremark{remark}{Remark}
\newsiamremark{hypothesis}{Hypothesis}
\crefname{hypothesis}{Hypothesis}{Hypotheses}
\newsiamthm{claim}{Claim}

\newtheorem{ass}[theorem]{Assumption}
\newtheorem{exm}[theorem]{Example}

\newcommand{\reff}[1]{(\ref{#1})}
\newcommand{\refa}[1]{Assumption\ \ref{#1}}
\newcommand{\refl}[1]{Lemma\ \ref{#1}}

\newcommand{\refd}[1]{Definition\ \ref{#1}}

\newcommand{\refal}[1]{Algorithm\ \ref{#1}}

\def\dd{&\,\,}
\def\st{\hbox{s.t.}}
\def\alglist{
	\begin{list}{Step 1}
		{\setlength{\leftmargin}{0.4 in}\setlength{\labelwidth}{0.7 in}}
	}
	\def\eli{\end{list}}

\def\na{\nabla}

\def\De{\Delta}

\def\la{\lambda}

\def\hf{\frac{1}{2}}

\def\alp{\alpha}

\title{Exact penalty method for D-stationary point of nonlinear optimization
}
\author{Xin-Wei Liu\thanks{Institute of Mathematics, Hebei University of Technology, Tianjin 300401, China
		(\email{mathlxw@hebut.edu.cn}). The research is supported by the NSFC grants (nos. 12071108 and 11671116).}
\and Yu-Hong Dai\thanks{Corresponding author. LSEC, Academy of Mathematics and Systems Science, Chinese Academy of Sciences, Beijing 100190, China; School of Mathematical Sciences, University of Chinese Academy of Sciences, Beijing 100049, China
  (\url{http://lsec.cc.ac.cn/\string~dyh/}). This author is supported by the NSFC grants (nos. 12021001, 11991021, 11991020 and 11971372) and  the Strategic Priority Research Program of Chinese Academy of Sciences (no. XDA27000000).}
}
\begin{document}

\maketitle

\begin{abstract}We consider the nonlinear optimization problem with least $\ell_1$-norm measure of constraint violations and introduce the concepts of the D-stationary point, the DL-stationary point and the DZ-stationary point with the help of exact penalty function.
If the stationary point is feasible, they correspond to the Fritz-John stationary point, the KKT stationary point and the singular stationary point, respectively. In order to show the usefulness of the new stationary points, we propose a new exact penalty sequential quadratic programming (SQP) method with inner and outer iterations and analyze its global and local convergence. The proposed method admits convergence to a D-stationary point and rapid infeasibility detection without driving the penalty parameter to zero, which demonstrates the commentary given in [SIAM J. Optim., 20 (2010), 2281--2299] and can be thought to be a supplement of the theory of nonlinear optimization on rapid detection of infeasibility. Some illustrative examples and preliminary numerical results demonstrate that the proposed method is robust and efficient in solving infeasible nonlinear problems and a degenerate problem without LICQ in the literature.
\end{abstract}



\pagestyle{myheadings}
\thispagestyle{plain}
\markboth{X.-W. LIU and Y.-H. DAI}
{Exact penalty method for D-stationary point of nonlinear optimization}

\section{Introduction}
\label{intro}

Consider the general nonlinear optimization \begin{eqnarray}
\min~ f(x) \quad\st~ h(x)=0,~g(x)\ge 0, \label{prob1-1}\end{eqnarray}
where $x\in\Re^n$ is the unknown, $f: \Re^n\to\Re$, $h: \Re^n\to\Re^{m_{\cal E}}$, $g: \Re^n\to\Re^{m_{\cal  I}}$ are twice continuously differentiable real-valued functions defined on $\Re^n$, $m_{\cal E}$ and $m_{\cal I}$ are two nonnegative integers with $\max\,\{m_{\cal E}, m_{\cal I}\}>0$. When $m_{\cal E}=0$ and $m_{\cal I}>0$, the problem is an inequality constrained optimization. If $m_{\cal E}>0$ and $m_{\cal I}=0$, the problem is an optimization problem with only equality constraints.

For the general problem \reff{prob1-1}, a {point} $x^*$ is called a {\bf Fritz-John point} if $x^*$ is feasible and there exist a scalar $\rho^*\ge 0$ and multiplier vectors $(\mu^*,\la^*)\in\Re^{m_{\cal E}}\times\Re_+^{m_{\cal I}}$ such that \begin{eqnarray}
\dd\dd \rho^*\na f(x^*)-\na h(x^*)\mu^*-\na g(x^*)\lambda^*=0, \label{wes1b}\\
\dd\dd \mu_i^*h_i(x^*)=0,~i=1,\ldots,m_{\cal E}, \\
\dd\dd \lambda_i^*g_i(x^*)=0,~i=1,\ldots,m_{\cal I}.\label{wes3b} \end{eqnarray}
Furthermore, if $\rho^*>0$, $x^*$ is called a {\bf KKT point}; if $\rho^*=0$, $x^*$ is called a {\bf singular stationary point}. A minimizer of problem \reff{prob1-1} is always a Fritz-John point, and is a KKT point under either the linear independence constraint qualification (LICQ) or the Mangasarian-Fromowitz constraint qualification (MFCQ).

Now let us consider the following simple one-dimensional example
 \begin{eqnarray} \min~x\quad\st~x^2-1\ge 0,~x-2\ge 0. \label{s-example}\end{eqnarray}
Surprisingly, when starting from a standard initial point $x_0=-4$, many methods including some state-of-the-art solvers can only converge to an infeasible point $x^*=-1$, as observed in \cite{ByrMaN,LiuDHS,LiuSun01,WacBie04}. How can we understand this point?

Since the point is even not feasible to problem \reff{s-example}, it is not a Fritz-John point at all. This means that, to have a comprehensive understandings on general nonlinear optimization problems, we have to face the fundamental feasibility/infeasibility issue. On one hand, the linearization of feasible nonlinear constraints may lead to inconsistent subproblems during the solution procedure (for example, see \cite{ByrCuN10,Char1978,ColCon80,Fletch87,GouOrT15,HanMan79,Piet1969}). On the other hand, the original nonlinear optimization problem may be infeasible in itself (for example, see \cite{BurCuW14,DLS17,LDHmm,NOW}).

Remarkably, to deal with the infeasibility issue, Burke and Han \cite{BurHan89} proposed the concept of {\bf infeasible stationary point}, which minimizes the constraint violation under some kind of measure. As Byrd, Curtis and Nocedal \cite{ByrCuN10} pointed out, there are various types of infeasible stationary points, some of which are strict isolated local minimizers of the infeasibility measure and others may belong to a set of minimizers of the measure or may simply be stationary points of the measure. Furthermore, Byrd, Curtis and Nocedal \cite{ByrCuN10} presented a set of conditions to guarantee the superlinear convergence of their SQP algorithm to an infeasible stationary point. Burke, Curtis and Wang \cite{BurCuW14} considered the general program with equality and inequality constraints, and proved that their SQP algorithm has strong global convergence and rapid convergence to a KKT point, and has superlinear/quadratic convergence to an infeasible stationary point. Recently, Dai, Liu and Sun \cite{DLS17} proposed a primal-dual interior-point method, which can be superlinearly or quadratically convergent to a KKT point if the original
problem is feasible, and can be superlinearly or quadratically convergent to an infeasible stationary point when the problem is infeasible.

As an infeasible stationary point is just to minimize the constraint violation under some kind of measure, it has nothing to do with the objective function. Motivated by the trajectory optimization problem in aerospace engineering,  Dai and Zhang \cite{DZ} proposed the basic model of optimization with least constraint violation. If the least $\ell_2$-norm measure of constraint violation is used, they provided
the optimality conditions for $M$-stationary point and $L$-stationary point and proved that the classical penalty method and the smoothing Fischer-Burmeister function method still work when the problem is infeasible. Marvelously, for convex optimization problem with least
constraint violation, Dai and Zhang \cite{DZa} demonstrated that the dual of the problem has an unbounded solution set, the optimality condition
can be characterized by the augmented Lagrangian, and the augmented Lagrangian method is able to find an approximate solution when the least violated shift is in the domain of the subdifferential of the optimal value function. More related work can be found in \cite{CD,ChiGil}.

Following the above line, we consider the nonlinear optimization with least $\ell_1$-norm measure of constraint violations. With the help of
exact penalty function (see Section 2 for details), we propose the following concept of $D$-stationary points for problem \reff{prob1-1}, which can be regarded as the generalizations of the well-known Fritz-John points, KKT points and singular stationary points of feasible optimization to infeasible optimization. A point $x^*$, which may be either feasible or infeasible to problem \reff{prob1-1}, is a called {\bf D-stationary point} of problem \reff{prob1-1} if it is a stationary point for minimizing the $\ell_1$-norm of constraint violations and there exist a scalar $\rho^*\ge 0$ and multiplier vectors $(\mu^*,\la^*)\in\Re^{m_{\cal E}}\times\Re_+^{m_{\cal I}}$ such that  \begin{eqnarray}
\dd\dd \rho^*\na f(x^*)-\na h(x^*)\mu^*-\na g(x^*)\lambda^*=0, \label{wes1a}\\
\dd\dd \mu_i^*h_i(x^*)+|h_i(x^*)|=0,~i=1,\ldots,m_{\cal E}, \\
\dd\dd \lambda_i^*g_i(x^*)+\max\{0,-g_i(x^*)\}=0,~i=1,\ldots,m_{\cal I}.\label{wes3a} \end{eqnarray}
Furthermore, if $\rho^*>0$, $x^*$ is called a {\bf DL-stationary point}. If $\rho^*=0$, $x^*$ is called a {\bf DZ-stationary point}.

The definitions of the D-stationary point, the DL-stationary point and the DZ-stationary point can be applied no matter whether problem \reff{prob1-1} is feasible or infeasible.
More exactly, we can see that a feasible D-stationary point is a Fritz-John point, a feasible DL-stationary point is a KKT point, and a feasible DZ-stationary point is a singular stationary point. In addition, an infeasible DZ-stationary point is an infeasible stationary point whose characteristics have been given in Burke and Han \cite{BurHan89} and Byrd, Curtis and Nocedal \cite{ByrCuN10}. For the simple one-dimensional problem \reff{s-example}, we can check that $x^*=-1$ is an infeasible D-stationary point (or an infeasible DL-stationary point)
with $\rho^*=0.1$ and $\lambda^*=(0.45, 1.0)$. Significantly, equations \reff{wes1a}--\reff{wes3a} correspond to the Fritz-John conditions of the feasible optimization with relaxation of least constraint violations (see \reff{relax5}).

In order to show the usefulness of our new stationary points, we shall propose a new exact penalty sequential quadratic programming (SQP) method with inner and outer iterations and analyze its global and local convergence. The proposed method generates inner iterations by an SQP method for the exact penalty subproblem with fixed penalty parameter and the outer iterations are responsible for update of the parameter. As a result, the method admits convergence to a D-stationary point and rapid infeasibility detection even though the penalty parameter is not small enough, which is distinguished from the existing methods with rapid detection of infeasibility such as \cite{BurCuW14,ByrCuN10,DLS17,LiuDHS}. The latter requires the penalty parameter tend to zero (or infinity as the sequence of parameters is nondecreasing). Some illustrative examples and preliminary numerical results demonstrate that the proposed method is robust and efficient in solving infeasible nonlinear problems and a degenerate problem without LICQ in the literature.

The paper is organized as follows. In section 2, the D-stationary point, DL-stationary point and DZ-stationary are introduced based on optimization with least $\ell_1$-norm constraint violation. Some simple illustrative examples are also given to show their relations to the exact penalty optimization. In section 3, we reformulate the original optimization problem into a new feasible problem with a parameter which is an smooth equivalent of the $\ell_1$ exact penalty optimization. Some basic results on the feasible reformulation and the optimization with least $\ell_1$-norm constraint violation are proved as well. We present our exact penalty SQP method in section 4. The global and local convergence are analyzed in sections 5 and 6, respectively. Some numerical results are reported in section 7. We conclude our paper in the last section.
	
Throughout the paper, we use standard notations from the literature. All vectors are column vectors. For $u, v\in\Re^s$, $\max\{u,v\}$ and $|u|$ are two $s$-dimensional vectors with components $\max\{u_i,v_i\}$ and $|u_i|$ for $i=1, \ldots, s$, respectively. For any vectors $x$ and $y$, $w=(x,y)$ means $w=[x^T\quad y^T]^T$. A letter with subscript $k$ is related to the $k$th iteration, the subscript $i$ indicates the $i$th component of a vector, and the subscript $kj$ is the $j$th iteration in solving the $k$th subproblem. The equation $\theta_k={O}(t_k)$ means that there exists a scalar $M$ independent of $k$ such that $|\theta_k|\le M|t_k|$ for all $k$ large enough, and equation $\theta_k={o}(t_k)$ means that $|\theta_k|\le\epsilon_k|t_k|$ for all $k$ large enough with $\lim_{k\to 0}\epsilon_k=0$. If it is not specified, letter $I$ is the identity matrix, ${\cal E}=\{1,2,\ldots,m_{\cal E}\}$ and ${\cal I}=\{1,2,\ldots,m_{\cal I}\}$ are two index sets, $\|\cdot\|$ be the Euclidean norm. Some unspecified notations may be identified from the context.

\section{D-stationary points for problem \reff{prob1-1} and some illustrative examples}

The stationary properties play a key role in recognizing optimal solutions and developing efficient algorithms for nonlinear optimization \reff{prob1-1}. For feasible optimization, the Fritz-John conditions and the KKT conditions are those of our well-known and extensively used stationary properties, where the Fritz-John conditions see \reff{wes1b}--\reff{wes3b}. The others are the conditions on singular stationary points. We restate the conditions on the KKT point and the singular stationary point so that the readers can have a better understanding on our introduction of new stationary points.
\begin{definition}\label{def0} Consider nonlinear optimization \reff{prob1-1}. Suppose that functions $f$, $g$ and $h$ are continuously differentiable on $\Re^n$, $x^*\in\Re^n$ is given. \\
(1) The point $x^*$ is a KKT point if it is feasible to the problem and there exist $\mu^*\in\Re^{m_{\cal E}}$ and $\la^*\in\Re_+^{m_{\cal I}}$ such that \begin{eqnarray}
\dd\dd \na f(x^*)-\na h(x^*)\mu^*-\na g(x^*)\lambda^*=0, \nonumber\\
\dd\dd \lambda_i^*g_i(x^*)=0,~i=1,\ldots,m_{\cal I}. \nonumber\end{eqnarray}
(2) The point $x^*$ is a singular stationary point if it is feasible to the problem and there exist $\mu^*\in\Re^{m_{\cal E}}$ and $\la^*\in\Re_+^{m_{\cal I}}$ such that $(\mu^*,\la^*)\ne 0$, and \begin{eqnarray}
\dd\dd -\na h(x^*)\mu^*-\na g(x^*)\lambda^*=0, \nonumber\\
\dd\dd \lambda_i^*g_i(x^*)=0,~i=1,\ldots,m_{\cal I}. \nonumber\end{eqnarray}\end{definition}

Based on the motivation to deal with infeasible optimization, we give the following definition on the DL-stationary point and the DZ-stationary point, which can be thought as a generalization of the conditions on the KKT point and the singular stationary point. In addition, some illustrative examples are given to show the significance of the stationary points for optimization with least $\ell_1$-norm constraint violations.
\begin{definition}\label{def1} Consider the nonlinear optimization problem \reff{prob1-1}. Suppose that functions $f$, $g$ and $h$ are continuously differentiable on $\Re^n$ and that $x^*\in\Re^n$ is a given point, which may be either feasible or infeasible to the problem. \\
(1) The point $x^*$ is referred to as a {\sl DL-stationary point} of nonlinear optimization \reff{prob1-1} if it is a stationary point of minimizing $\|h(x)\|_1+\|\max\{0, -g(x)\}\|_1$, and there exist
multiplers $(\mu^*,\lambda^*)\in\Re^{m_{\cal E}}\times\Re_+^{m_{\cal I}}$ such that \begin{eqnarray}
\dd\dd \na f(x^*)-\na h(x^*)\mu^*-\na g(x^*)\lambda^*=0, \label{wes1}\\
\dd\dd \mu_i^*h_i(x^*)=-|h_i(x^*)|,~i=1,\ldots,m_{\cal E}, \\
\dd\dd \lambda_i^*g_i(x^*)=-\max\{0,-g_i(x^*)\},~i=1,\ldots,m_{\cal I}.\label{wes3} \end{eqnarray}
(2) The point $x^*$ is a {\sl DZ-stationary point} if it is a stationary point for minimizing $\|h(x)\|_1+\|\max\{0, -g(x)\}\|_1$, and there exist $\mu^*\in\Re^{m_{\cal E}}$ and $\la^*\in\Re_+^{m_{\cal I}}$ such that \begin{eqnarray}
\dd\dd -\na h(x^*)\mu^*-\na g(x^*)\lambda^*=0, \label{wes1c}\\
\dd\dd \mu_i^*h_i(x^*)=-|h_i(x^*)|,~i=1,\ldots,m_{\cal E}, \\ 
\dd\dd \lambda_i^*g_i(x^*)=-\max\{0,-g_i(x^*)\},~i=1,\ldots,m_{\cal I}. \label{wes3c}\end{eqnarray}
\end{definition}

Similar to that both the KKT point and the singular stationary point are the Fritz-John point, both the DL-stationary point and the DZ-stationary point belong to the set of D-stationary points. In essence, if $x^*$ is a stationary point for minimizing $\|h(x)\|_1+\|\max\{0, -g(x)\}\|_1$, then the conditions \reff{wes1c}--\reff{wes3c} hold naturally with some $\mu^*\in\Re^{m_{\cal E}}$ and $\la^*\in\Re_+^{m_{\cal I}}$. Thus, a  DZ-stationary point is a singular stationary point if the point is feasible and it is an infeasible stationary point if the point is infeasible.
What is more interesting and important, when the optimization problem has inherent violations $(v_1^h,\ldots,v_{m_{\cal E}}^h)$ for equality constraints and $(v_1^g,\ldots,v_{m_{\cal I}}^g)$ for inequality constraints, the conditions \reff{wes1}--\reff{wes3} with all $|h_i(x^*)|$ and $\max\{0,-g_i(x^*)\}$ being replaced, respectively, by $v_i^h$ and $v_i^g$ correspond to the KKT conditions of the feasible optimization, a relaxation of the original problem,  \begin{eqnarray}\min~f(x)\quad
\st~h_i(x)=v_i^h,~\forall i\in{\cal E};~g_i(x)\ge v_i^g,~\forall i\in{\cal I}. \label{relax5}\end{eqnarray}
Thus we can see that D-stationary point, DL-stationary point and DZ-stationary point are generalization of Fritz-John point, KKT point and singular stationary point, respectively, from the feasible case to the general case.

In order to show the usefulness of our new stationary points, based on robustness of the exact penalty optimization, we investigate the stationary properties of an exact penalty optimization. To this aim, consider the $\ell_1$ exact penalty approach on the original nonlinear optimization \reff{prob1-1},
\begin{eqnarray} \min_{x\in\Re^n}~P_{\rho}(x)=:\rho f(x)+c(x), \label{prob1L1}\end{eqnarray}
where $c: \Re^n\to\Re$ and $c(x)=\|h(x)\|_1+\|\max\{0, -g(x)\}\|_1$ is the $\ell_1$-norm measure of constraint violations,  $\rho>0$ is a penalty parameter. Apparently, $c(x)$ is nonsmooth even though both $h(x)$ and $g(x)$ are twice differentiable. Moreover, $c(x)\ge 0$ for all $x\in\Re^n$. It is well-known that $P_{\rho}(x)$ is an exact penalty function in the sense that,
if the original problem \reff{prob1-1} is feasible, its local minima can be the minima of problem \reff{prob1L1} under second-order sufficient conditions for all penalty parameters smaller than certain threshold.

We further reformulate the non-smooth $\ell_1$ exact penalty optimization \reff{prob1L1} as a smooth feasible problem with only inequality constraints,  \begin{eqnarray}
\min_{x,y,z}\dd\dd\rho f(x)+\sum_{i=1}^{m_{\cal E}} y_i+\sum_{i=1}^{m_{\cal I}} z_i, \label{prob2-1}\\
\st\dd\dd y-h(x)\ge 0,~ y+h(x)\ge 0, \label{prob2-2}\label{prob2-3}\\
   \dd\dd z+g(x)\ge 0,~ z\ge 0, \label{prob2-4}\label{prob2-5}\end{eqnarray}
where $x\in\Re^n$, $y\in\Re^{m_{\cal E}}$, $z\in\Re^{m_{\cal I}}$, $\rho>0$ is the penalty parameter. Since the constraints \reff{prob2-3} and \reff{prob2-5} guarantee $y\ge 0$ and $z\ge 0$, one has $\|y\|_1=\sum_{i=1}^{m_{\cal E}}y_i$ and $\|z\|_1=\sum_{i=1}^{m_{\cal I}}z_i$. The same reformulation appears in (2.5) of \cite{GouOrT15} and is an extension of (2.2) in \cite{ByrCuN10}. It is noted that \cite{BurCuW14,NocWri99} use a different reformulation with both equality and inequality constraints for $\ell_1$ exact penalty problem \reff{prob1L1}, where they need double not single additional nonnegative variable as ours for every equality constraint.
For convenience of statement, we denote \begin{eqnarray}
\dd\dd F_{\rho}(x,y,z)=\rho f(x)+\sum_{i=1}^{m_{\cal E}} y_i+\sum_{i=1}^{m_{\cal I}} z_i, \nonumber\\
\dd\dd G(x,y,z)=(y-h(x),y+h(x),z+g(x),z). \nonumber\end{eqnarray}

Let us look at three simple examples, by which we may have an intuitive understanding on our new stationary points and their relations to the original optimization, the exact penalty optimization and the optimization with least $\ell_1$ constraint violations.

\begin{exm}\label{ex1} Consider the nonlinear constrained optimization problem \begin{eqnarray}
\min_x~f(x):=x^2+4x \quad\st~ h(x):=x-1=0. \nonumber\end{eqnarray}\end{exm}

Obviously, the problem is feasible and has a solution $x^*=1$. The reformulation \reff{prob2-1}--\reff{prob2-5} has the form \begin{eqnarray}
\min_{x,y}\dd\dd F_{\rho}(x,y):=\rho x^2+4\rho x+y \nonumber\\
\st\dd\dd G(x,y):=\left[\begin{array}{l}
y-x+1 \\
y+x-1 \end{array}\right]\ge 0. \nonumber\label{probe2-2}\end{eqnarray}
If $\rho>\frac{1}{6}$, one has $(x^*,y^*)=(\frac{1}{2\rho}-2,3-\frac{1}{2\rho})$, $(u^*,v^*)=(0,1)$;
otherwise, $(x^*,y^*)=(1,0)$, $(u^*,v^*)=(\frac{1}{2}-3\rho,\frac{1}{2}+3\rho)$. Moreover, \begin{eqnarray}
F_{\rho}(x^*,y^*)=\left\{\begin{array}{ll}
-\frac{1}{4\rho}-4\rho+3, & \rho>\frac{1}{6}; \\[2mm]
5\rho, & \rho\le\frac{1}{6}. \end{array}\right. \nonumber\end{eqnarray}
Finally, $x^*=1$ is a KKT point of the problem. The associated Lagrange multiplier is $\mu^*=6$. The threshold is $\rho^*=\frac{1}{6}$.

\begin{exm} Consider an infeasible nonlinear constrained optimization problem \begin{eqnarray}
\min_x~f(x):=x^2+4x\quad\st~ h(x):=\left[\begin{array}{l}
x-1\\
x+1 \end{array}\right]=0. \nonumber\end{eqnarray}\end{exm}

Corresponding to the reformulation \reff{prob2-1}--\reff{prob2-5}, we have the problem \begin{eqnarray}
\min_{x,y}\dd\dd F_{\rho}(x,y):=\rho x^2+4\rho x+(y_1+y_2) \nonumber\\
\st\dd\dd G(x,y):=\left[\begin{array}{l}
y_1-x+1 \\
y_1+x-1 \\
y_2-x-1 \\
y_2+x+1 \end{array}\right]\ge 0. \nonumber\end{eqnarray}
If $\rho>1$, the KKT point of the reformulation problem is $(x^*,y_1^*,y_2^*)=(\frac{1}{\rho}-2,3-\frac{1}{\rho},1-\frac{1}{\rho})$, and $(u_1^*,v_1^*,u_2^*,v_2^*)=(0,1,0,1)$ is the associated Lagrange multiplier vector.
If $\rho\le1$, the KKT point and the associated Lagrange multiplier are $(x^*,y_1^*,y_2^*)=(-1, 2, 0)$ and $(u_1^*,v_1^*,u_2^*,v_2^*)=(0,1,1-\rho,\rho)$, respectively.
Thus, \begin{eqnarray}
F_{\rho}(x^*,y^*)=\left\{\begin{array}{ll}
-\frac{1}{\rho}-4\rho+4, & \rho>1; \\[2mm]
2-3\rho, & \rho\le1. \end{array}\right. \nonumber\end{eqnarray}

Note that $[-1,1]$ is the set of solutions for minimizing $|x-1|+|x+1|$. By \refd{def1}, both $x^*=-1$ and $x^*=1$ are the DL-stationary points, where our solution $x^*=-1$ is the unique global minimizer of optimization with least $\ell_1$-norm constraint violations \begin{eqnarray}
\min~ x^2+4x\quad\hbox{s.t.}~ x~\hbox{minimizes}~|x-1|+|x+1|. \nonumber\end{eqnarray}

\begin{exm} Consider an inequality constrained optimization with inconsistent constraints (see also Dai and Zhang \cite{DZ}) \begin{eqnarray}
\min~ x_1^2+x_2^2 \quad\st~ -x_1-x_2+1\ge 0, ~ x_1+x_2-2\ge 0. \nonumber\end{eqnarray}\end{exm}

The exact penalty approach solves the problem by solving its smooth and feasible reformulation \begin{eqnarray}
\min\dd\dd \rho x_1^2+\rho x_2^2+z_1+z_2 \nonumber\\
\st\dd\dd z_1-x_1-x_2+1\ge 0, \nonumber\\
\dd\dd z_2+x_1+x_2-2\ge 0, \nonumber\\
\dd\dd z_1\ge 0,~z_2\ge 0. \nonumber\end{eqnarray}
Although the original problem is infeasible, the preceding problem is well posed. By solving the KKT conditions, we have $x_1^*=x_2^*=0.5$, $z_1^*=0$, $z_2^*=1$, $s_1^*=1-\rho$, $s_2^*=1$, $t_1^*=\rho$, $t_2^*=0$ provided $\rho\le 1$. This means that, for $0<\rho\le 1$, we find a point with $f(x^*)=0.5$ at which the $\ell_1$-norm and the $\ell_2$-norm of constraint violations are $1$.

In contrast, the method in \cite{DZ} got the solution $x_1^*=x_2^*=0.75$, $f(x^*)=\frac{9}{8}$, which has the same $\ell_1$-norm constraint violations with ours and has less $\ell_2$-norm constraint violations $\frac{\sqrt{2}}{2}$ and larger objective value $\frac{9}{8}$. In other words, our method finds a point with much better objective value than the one proposed in \cite{DZ} under the same $\ell_1$-norm constraint violations. In addition, if $\rho>1$, one has $x^*=(\frac{1}{2\rho},\frac{1}{2\rho})$, $z^*=(0,2-\frac{1}{\rho})$, $s^*=(0,1)$, $t^*=(1,0)$. Finally, one has \begin{eqnarray}
F_{\rho}(x^*,z^*)=\left\{\begin{array}{ll}
-\frac{1}{2\rho}+2, & \rho>1; \\[2mm]
0.5\rho+1, & \rho\le 1. \end{array}\right. \nonumber\end{eqnarray}
It is easy to examine that all points satisfying one of the two constraints are the ones minimizing the $\ell_1$-norm constraint violations, among which our solution $(0.5,0.5)$ is a DL-stationary point and is the unique global minimizer of optimization with least $\ell_1$-norm constraint violations.

\section{The exact penalty optimization and the optimization with least constraint violations}

In this section, we investigate the relations among the original nonlinear optimization, the exact penalty optimization and its smooth reformulation, and the optimization with least constraint violations theoretically. Moreover, we provide a certificate on the solution of optimization with least $\ell_1$-norm constraint violations. The certificate is important in addressing convergence to the D-stationary points and the rapid infeasibility detection as the penalty parameter is not small.

Corresponding to the problem \reff{prob1-1}, the optimization with least $\ell_1$-norm measure of constraint violations is a bi-level optimization problem \begin{eqnarray}
\min~ f(x)\quad\st~ x\in\hbox{argmin}_{x\in\Re^n}~c(x):=\|h(x)\|_1+\|\max\{0, -g(x)\}\|_1. \label{prob1b1-2}\label{prob1b1-1}\end{eqnarray}
If problem \reff{prob1-1} is feasible, then it is equivalent to the problem \reff{prob1b1-1} in sense that both of them have the same minima. However, like problem \reff{prob1L1} but not the original problem, the formulation \reff{prob1b1-1} admits infeasible constraints.
Since it is subject to the set of solutions of an unconstrained non-smooth problem, the formulation \reff{prob1b1-1} is difficult and complicated in how to describe its optimality conditions and design efficient algorithms. These difficulties may be alleviated when the original problem \reff{prob1-1} is convex, see \cite{CD,ChiGil,DZ,DZa} and references therein. For general nonconvex nonlinear optimization, we investigate the stationary points of problem \reff{prob1b1-1} with the help of exact penalty optimization \reff{prob1L1}, which results in our new stationary points of nonlinear optimization.

The following results show the relations between exact penalty optimization \reff{prob1L1} and the original problem, where the second part of results can see Theorem 4.4 of \cite{HanMan79} and Theorem 17.3 of \cite{NocWri99}.
\begin{proposition} If $x^*$ is a local minimizer of problem \reff{prob1L1} and $c(x^*)=0$, $\rho>0$ is any penalty parameter, then $x^*$ is a local minimizer of the original nonlinear optimization \reff{prob1-1}. Conversely, if $x^*$ is a strict local minimizer of the original problem \reff{prob1-1} at which the KKT conditions are satisfied with Lagrange multipliers $\mu^*\in\Re^{m_{\cal E}}$ and $\la^*\in \Re_+^{m_{\cal I}}$, then $x^*$ is a local minimizer of problem \reff{prob1L1} for all $\rho<\rho^*$, where $\rho^*=\frac{1}{\max\{\|\mu^*\|_{\infty},\|\la\|^*_{\infty}\}}$.
If, in addition, the second-order sufficient conditions of problem \reff{prob1-1} hold at the triple $(x^*,\mu^*,\la^*)$ and $\rho<\rho^*$, then $x^*$ is a strict local minimizer of problem \reff{prob1L1}. \end{proposition}\begin{proof}
If the point $x^*$ is a local minimizer of problem \reff{prob1L1}, we have for some $\delta>0$ and all $x\in{\cal N}(x^*,\delta)$, \begin{eqnarray}
\rho f(x^*)+c(x^*)\le\rho f(x)+c(x), \nonumber\end{eqnarray}
where $\delta$ is a scalar and ${\cal N}(x^*,\delta)$ is a $\delta$-neighborhood of $x^*$. If $c(x^*)=0$, the preceding inequality is reduced to \begin{eqnarray}
\rho f(x^*)\le\rho f(x)+c(x),\quad\forall x\in{\cal N}(x^*,\delta). \nonumber\end{eqnarray}
Hence, $f(x^*)\le f(x)$ for those $x\in{\cal N}(x^*,\delta)\cap\{x|c(x)=0\}$. That is, $x^*$ is a local minimizer of nonlinear optimization \reff{prob1-1}.

The proof for the remaining results can see that of Theorem 4.4 of \cite{HanMan79}.
\end{proof}

Let us revisit problem \reff{prob2-1}--\reff{prob2-5}, the smooth equivalent of problem \reff{prob1L1}. Apparently, for every $x\in\Re^n$, $(x,y,z)$ with $y=|h(x)|$ and $z=\max\{0,-g(x)\}$ is a feasible point of problem \reff{prob2-1}--\reff{prob2-5}. Thus, problem \reff{prob2-1}--\reff{prob2-5} is always feasible no matter whether problem \reff{prob1-1} is feasible or not. Similar to the reformulation (1.3) of \cite{GouOrT15}, we have the regular properties on problem \reff{prob2-1}--\reff{prob2-5} as follows.
\begin{proposition}\label{lemma1} Problem \reff{prob2-1}--\reff{prob2-5} is always feasible. Moreover, the following results are true. \\
(1) MFCQ holds at any feasible point $(x,y,z)$ of problem \reff{prob2-1}--\reff{prob2-5}; \\
(2) If a point $\hat x$ is feasible to the problem \reff{prob1-1} and LICQ holds at $\hat x$, then $(\hat x,\hat y,\hat z)$ with $\hat y=0$ and $\hat z=0$ is a feasible point of problem \reff{prob2-1}--\reff{prob2-5} and LICQ holds at $(\hat x,\hat y,\hat z)$; \\
(3) Let $(\hat x,\hat y,\hat z)$ be a feasible point of problem \reff{prob2-1}--\reff{prob2-5}, $\hat{\cal I}_h=\{i| h_i(\hat x)=0\}$, $\hat{\cal I}_g=\{i| g_i(\hat x)=0\}$. If $\na h_i(\hat x) (i\in\hat{\cal I}_h)$, $\na g_i(\hat x) (i\in\hat{\cal I}_g)$ are linearly independent, then for problem \reff{prob2-1}--\reff{prob2-5}, LICQ holds at $(\hat x,\hat y,\hat z)$.
\end{proposition}\begin{proof} For every $x\in\Re^n$, let $y=|h(x)|$, $z=\max\{0, -g(x)\}$. Then all constraints in problem \reff{prob2-1}--\reff{prob2-5} are satisfied with $(x,y,z)$. Thus, problem \reff{prob2-1}--\reff{prob2-5} is feasible for all $x\in\Re^n$ and suitably selected $y\in\Re^{m_{\cal E}}$ and $z\in\Re^{m_{\cal I}}$.

(1) Note that \begin{eqnarray} \na G(x,y,z)=\left[\begin{array}{cccc}
-\na h(x) & \na h(x) & \na g(x) & 0 \\
I_{m_{\cal E}} & I_{m_{\cal E}} & 0 & 0 \\
0 & 0 & I_{m_{\cal I}} & I_{m_{\cal I}}\end{array}\right]. \nonumber\end{eqnarray}
For any $\De w=(\De x,\De y,\De z)$, \begin{eqnarray} \na G(x,y,z)^T\De w>0~\hbox{if and only if}~\left\{\begin{array}{l}
\De y-\na h(x)^T\De x>0,  \\
\De y+\na h(x)^T\De x>0,  \\
\De z+\na g(x)^T\De x>0,  \\
\De z>0. \end{array}\right.\nonumber\end{eqnarray}
Thus, for any $\De x$, there are $\De y>|\na h(x)^T\De x|$ and $\De z>\max\{0,-\na g(x)^T\De x\}$ such that $\na G(x,y,z)^T\De w>0$, which implies that MFCQ holds immediately.

(2) and (3) can be verified by examining the columns of $\na G(x,y,z)$.  \end{proof}

Problem \reff{prob2-1}--\reff{prob2-5} is robust in the sense that it is feasible even though the constraints in the original problem \reff{prob1-1} are inconsistent, and its local minima are always the KKT points since MFCQ naturally holds at any feasible points. Furthermore, the associated Lagrange multipliers of any KKT points of the smooth problem \reff{prob2-1}--\reff{prob2-5} are bounded. \begin{proposition}\label{pro230116a}
If the exact penalty funtion $P_{\rho}:\Re^n\to\Re$ is lower bounded and $x^*$ is a local minimizer, then there exist Lagrange multipliers $u^*\in\Re_+^{m_{\cal E}}$, $v^*\in\Re_+^{m_{\cal E}}$, $s^*\in\Re_+^{m_{\cal I}}$, $t^*\in\Re_+^{m_{\cal I}}$ such that the KKT conditions of problem \reff{prob2-1}--\reff{prob2-5} hold at $(x^*,y^*,z^*)$ with $y^*=|h(x^*)|$ and $z^*=\max\{0,-g(x^*)\}$. Moreover, $u^*\le e_{m_{\cal E}}$, $v^*\le e_{m_{\cal E}}$, $s^*\le e_{m_{\cal I}}$, $t^*\le e_{m_{\cal I}}$, where $e_{m_{\cal E}}\in\Re^{m_{\cal E}}$ and $e_{m_{\cal I}}\in\Re^{m_{\cal I}}$ are all-one vectors. \end{proposition}\begin{proof}
Due to the equivalence between problems \reff{prob1L1} and \reff{prob2-1}--\reff{prob2-5}, $x^*$ is a minimizer of problem \reff{prob1L1} if and only if $(x^*,y^*,z^*)$ is a minimizer of problem \reff{prob2-1}--\reff{prob2-5}.
Since MFCQ holds at $(x^*,y^*,z^*)$, it follows from the first-order necessary optimality conditions of general constrained optimization (see, for example, \cite{NocWri99,SunYua06}) that any minimizer $(x^*,y^*,z^*)$ of problem \reff{prob2-1}--\reff{prob2-5} satisfies the KKT conditions \begin{eqnarray}
\dd\dd\rho\na f(x^*)+\na h(x^*)(u^*-v^*)-\na g(x^*)s^*=0, \label{kkt2.9-1}\\
\dd\dd 1-u_i^*-v_i^*=0,\ i=1,\ldots,m_{\cal E}, \label{kkt2.9-2}\\
\dd\dd 1-s_i^*-t_i^*=0,\ i=1,\ldots,m_{\cal I}, \label{kkt2.9-3}\\
\dd\dd (u^*)^T(y^*-h(x^*))=0,~ (v^*)^T(y^*+h(x^*))=0, \\
\dd\dd (s^*)^T(z^*+g(x^*))=0,~ (t^*)^Tz^*=0, \label{kkt2.9-8}\end{eqnarray}
where $u^*\ge 0$, $v^*\ge 0$, $s^*\ge 0$, $t^*\ge 0$ are associated Lagrange multipliers. The boundedness of $u^*$, $v^*$, $s^*$ and $t^*$ derive from equations \reff{kkt2.9-2}, \reff{kkt2.9-3} and the nonnegativity immediately. \end{proof}

Similar results have been presented in \cite{BurCuW14,ByrCuN10,GouOrT15}. For example, Gould et al. in \cite{GouOrT15} think that the reformulation is surprisingly regular.
Let $L(x,\mu,\lambda)=f(x)-\mu^Th(x)-\lambda^Tg(x)$ be the Lagrangian of problem \reff{prob1-1}, where $\mu$ and $\lambda$ are Lagrange multipliers associated with the equality and inequality constraints, respectively. In what follows, we describe the relations of minimizers between problem \reff{prob1-1} and its robust counterpart \reff{prob2-1}--\reff{prob2-5}.
\begin{proposition}\label{lemma2} Let $x^*$ be a local minimizer of problem \reff{prob1-1} at which LICQ and the second-order sufficient conditions (SOSCs) are satisfied with Lagrange multipliers $\mu^*\in\Re^{m_{\cal E}}$ and $\lambda^*\in\Re_{+}^{m_{\cal I}}$. Then there exists a threshold value $\rho^*>0$ such that, for all $0<\rho\le\rho^*$, $(x^*,y^*,z^*)$ together with $y^*=0$ and $z^*=0$ is a strict local minimizer of problem \reff{prob2-1}--\reff{prob2-5}, and there exist Lagrange multipliers $u^*\in\Re_+^{m_{\cal E}}$, $v^*\in\Re_+^{m_{\cal E}}$, $s^*\in\Re_+^{m_{\cal I}}$, $t^*\in\Re_+^{m_{\cal I}}$ such that for problem \reff{prob2-1}--\reff{prob2-5}, LICQ, the KKT conditions and SOSCs are satisfied at $(x^*,y^*,z^*)$.
\end{proposition}\begin{proof} The result follows from the equivalence of problem \reff{prob2-1}--\reff{prob2-5} to the exact penalty optimization and Theorem 17.3 of \cite{NocWri99}. In particular, $(x^*,0,0)$ is feasible to the problem \reff{prob2-1}--\reff{prob2-5}, and the LICQ at $(x^*,0,0)$ is implied by LICQ at $x^*$ for problem \reff{prob1-1}. The penalty parameter $\rho^*=\frac{1}{\max\{\|\mu^*\|_{\infty},\|\lambda^*\|_{\infty}\}}$, and
$u_i^*=\hf{(1-\rho\mu_i^*)},~v_i^*=\hf{(1+\rho\mu_i^*)}~\hbox{for}~ i=1,\ldots,m_{\cal E},$
$s_i^*={\rho}{\lambda_i^*},~t_i^*=1-{\rho\lambda_i^*}~\hbox{for}~i=1,\ldots,m_{\cal I}.$ \end{proof}

The preceding proposition shows that, when the original nonlinear constrained optimization \reff{prob1-1} is feasible, the reformulation problem \reff{prob2-1}--\reff{prob2-5} can be equivalent to it when the parameter $\rho \le\frac{1}{\max\{\|\mu^*\|_{\infty},\|\lambda^*\|_{\infty}\}}$. However, this threshold is unknown before solving the original problem, and certain adaptive update procedure on the penalty parameter should be introduced to promote convergence of the methods based on solving the reformulation problem.

The following result is important since it provides a certificate on that a solution of the reformulation \reff{prob2-1}--\reff{prob2-5} can be a solution of optimization with least $\ell_1$-norm constraint violations \reff{prob1b1-1}. If the solution is not a solution of an equality constrained optimization problem based the active set of constraints at the solution of problem \reff{prob2-1}--\reff{prob2-5}, it is not a solution of optimization with least $\ell_1$-norm constraint violations. Thus, it establishes a foundation for designing methods for nonlinear optimization based on optimization with least $\ell_1$-norm constraint violations.
\begin{proposition}\label{theorem2} For any fixed $\rho>0$, suppose that $(x^*,y^*,z^*)$ with $y^*=|h(x^*)|$ and $z^*=\max\{0, -g(x^*)\}$ is a KKT point of problem \reff{prob2-1}--\reff{prob2-5}, $(u^*,v^*,s^*,t^*)$ with $u^*\in\Re_+^{m_{\cal E}}$, $v^*\in\Re_+^{m_{\cal E}}$, $s^*\in\Re_+^{m_{\cal I}}$, $t^*\in\Re_+^{m_{\cal I}}$ are the associated Lagrange multipliers. Let ${\cal I}_h^*=\{i| h_i(x^*)=0\}$, ${\cal I}_g^*=\{i| g_i(x^*)=0\}$.
If $x^*$ is a local solution of problem
\begin{eqnarray} \min_{x\in\Re^n}~\|h(x)\|_1+\|\max\{0, -g(x)\}\|_1, \label{L1p}\end{eqnarray}
then it is a KKT point of problem \begin{eqnarray}
\min_x\dd\dd f(x) \label{kkt2.13a-1}\\
\st\dd\dd h_i(x)=0,~i\in{\cal I}_h^*, \label{kkt2.13a-2}\\
\dd\dd g_i(x)=0,~i\in{\cal I}_g^*. \label{kkt2.13a-3}\end{eqnarray}
Furthermore, if gradients $\na h_i(x^*) (i\in{\cal I}_h^*)$, $\na g_i(x^*) (i\in{\cal I}_g^*)$ are linearly independent, and there holds $d^T\na_{xx}^2\hat L(x^*,\hat\mu^*,\hat\la^*)d>0$ for all $d\ne 0$ satisfying $\na h_i(x^*)^Td=0$ for $i\in{\cal I}_h^*$ and $\na g_i(x^*)^Td=0$ for $i\in{\cal I}_g^*$, where $\hat\mu_i^*  (i\in{\cal I}_h^*)$ and $\hat\la_i^* (i\in{\cal I}_g^*)$ are the Lagrange multipliers of problem \reff{kkt2.13a-1}--\reff{kkt2.13a-3} and $\hat L(x^*,\hat\mu^*,\hat\la^*)=f(x^*)-\sum_{i\in{\cal I}_h^*}\hat\mu_i^* h_i(x^*)-\sum_{i\in{\cal I}_g^*}\hat\la_i^* g_i(x^*)$, then $x^*$ is a strict local minimizer of optimization with least $\ell_1$-norm constraint violations \reff{prob1b1-1} (i.e., $x^*$ is a local solution for minimizing $\ell_1$-norm constraint violations and is a strict local minimizer of $f(x)$ in certain neighborhood of $x^*$). \end{proposition}\begin{proof}
Note that $(x^*,y^*,z^*,u^*,v^*,s^*,t^*)$ satisfies the KKT conditions \reff{kkt2.9-1}--\reff{kkt2.9-8} of problem \reff{prob2-1}--\reff{prob2-5}.
Then $u_i^*=1,~v_i^*=0$ for $i\in\{i|h_i(x^*)>0\}$, $u_i^*=0,~v_i^*=1$ for $i\in\{i|h_i(x^*)<0\}$, $s_i^*=0,~t_i^*=1$ for $i\in\{i|g_i(x^*)>0\}$, $s_i^*=1,~t_i^*=0$ for $i\in\{i|g_i(x^*)<0\}$.

The fact that $x^*$ is a solution of problem \reff{L1p} implies that there exist $\xi_i^*\in [-1,1], i\in{\cal I}_h^*$ and $\eta_i^*\in [-1,0], i\in{\cal I}_g^*$ such that \begin{eqnarray}\dd\dd\sum_{i\in\{i|h_i(x^*)>0\}}\na h_i(x^*)-\sum_{i\in\{i|h_i(x^*)<0\}}\na h_i(x^*)+\sum_{i\in{\cal I}_h^*}\xi_i^*\na h_i(x^*) \nonumber\\
\dd\dd\quad\quad-\sum_{i\in\{i|g_i(x^*)<0\}}\na g_i(x^*)+\sum_{i\in{\cal I}_g^*}\eta_i^*\na g_i(x^*)=0. \label{210912a}\end{eqnarray}
Subtracting \reff{210912a} from the left-hand-side of \reff{kkt2.9-1}, we derive \begin{eqnarray}
\rho\na f(x^*)-\sum_{i\in{\cal I}_h^*}\tilde\mu_i^*\na h_i(x^*)-\sum_{i\in{\cal I}_g^*}\tilde\la_i^*\na g_i(x^*)=0, \nonumber\end{eqnarray}
where $\tilde\mu_i^*=\xi_i^*-u_i^*+v_i^*$, $\tilde\la_i^*=\eta_i^*+s_i^*$. Thus, $x^*$ is a KKT point of problem \reff{kkt2.13a-1}--\reff{kkt2.13a-3}, $\hat\mu_i^*=\frac{\tilde\mu_i^*}{\rho},  i\in{\cal I}_h^*$ and $\hat\la_i^*=\frac{\tilde\la_i^*}{\rho},  i\in{\cal I}_g^*$ are the associated Lagrange multipliers.

It is well known that conditions of the proposition will guarantee $x^*$ to be a strict local minimizer of problem \reff{kkt2.13a-1}--\reff{kkt2.13a-3} (see \cite{NocWri99,SunYua06}).
Since by continuity all points around $x^*$ is a subset of the feasible set of problem \reff{kkt2.13a-1}--\reff{kkt2.13a-3}, the desired result follows immediately. \end{proof}

\section{An exact penalty SQP algorithm}

The SQP approach is very effective in solving general nonlinear optimization. In this section, we present an exact penalty SQP method for the original problem \reff{prob1-1} based on solving the subproblem \reff{prob2-1}--\reff{prob2-5}. Our method can be taken as a modified line search S$\ell_1$QP which not only is applicable to the feasible nonlinear constrained optimization like the generic S$\ell_1$QP (see \cite{Fletch87}), but also is applicable to the infeasible problems.
It is also a modification of the method \cite{BurCuW14} where we use a simple rule for penalty update and inner-outer iterative framework to replace the sophisticated rule and single iterative framework. Moreover, our QP subproblems are different from those of \cite{BurCuW14} in coping with equality constraints.
The proposed method shares similarity to the methods in \cite{BurCuW14,ByrCuN10} which are based on solving QP subproblems and can be rapidly detecting the infeasibility of nonlinear constrained optimization. It is also similar to that of \cite{ByrCuN10} which may be necessary for solving several QP subproblems before updating the penalty parameter, although the aims and results are not the same where we use every direction generated by QP subproblems to obtain a new inner iteration point.
Our gains from these modifications are the global and local convergence of the exact penalty approach to the stationary points of nonlinear optimization with least constraint violations.

For problem \reff{prob2-1}--\reff{prob2-5}, general SQP solves the QP subproblem \begin{eqnarray}
\min_{p\in\Re^{n}\times\Re^{m_{\cal E}}\times\Re^{m_{\cal I}}}\dd\dd \na F_{\rho}(x,y,z)^Tp+\hf p^TH p \label{210706a}\\
\st\dd\dd G(x,y,z)+\na G(x,y,z)^Tp\ge 0, \label{210706b}\end{eqnarray}
where $H$ can be the exact Lagrangian Hessian $$\left[\begin{array}{cc}
H_{\rho}(x,u,v,s) & 0_{n\times ({m}_{\cal E}+m_{\cal I})} \\
0_{({m}_{\cal E}+m_{\cal I})\times n} & 0_{({m}_{\cal E}+m_{\cal I})\times ({m}_{\cal E}+m_{\cal I})} \end{array}\right]$$
with $H_{\rho}(x,u,v,s)=\rho\na^2 f(x)+\sum_{i=1}^{m_{\cal E}}(u_i-v_i)\na^2 h_i(x)-\sum_{i=1}^{m_{\cal I}}s_i\na^2 g_i(x)\in\Re^{n\times n}$ or its appropriate approximation, $(x,y,z)$ is the current iterate, $(u,v,s,t)$ is the estimate of the associated Lagrange multiplier vector. In particular, let $p=(d,y^+-y,z^+-z)$, $B_{\rho}$ be a positive definite approximation to $H_{\rho}(x,u,v,s)$, then subproblem \reff{210706a}--\reff{210706b} is the following feasible QP \begin{eqnarray}
\min_{d,y^+,z^+}\dd\dd e_{m_{\cal E}}^Ty^++e_{m_{\cal I}}^Tz^++\rho \na f(x)^Td+\hf d^TB_{\rho}d \label{qp1a1}\\
\st\dd\dd y^+-h(x)-\na h(x)^Td\ge 0, \label{qp1a2}\\
\dd\dd y^++h(x)+\na h(x)^Td\ge 0, \label{qp1a3}\\
\dd\dd z^++g(x)+\na g(x)^Td\ge 0, \label{qp1a4}\\
\dd\dd z^+\ge 0, \label{qp1a5}\end{eqnarray}
where $e_{m_{\cal E}}$ and $e_{m_{\cal I}}$ are all-one vectors. Let $$c'(x;d)=\|h(x)+\na h(x)^Td\|_1+\|\max\{0,-(g(x)+\na g(x)^Td)\}\|_1.$$
The preceding problem is an equivalent of non-smooth $\ell_1$QP subproblem of \cite{Fletch87}, \begin{eqnarray}
\min_d\dd\dd \rho\na f(x)^Td+\hf d^TB_{\rho}d+c'(x;d). \nonumber\end{eqnarray}

If $(d_{\ell},y_{\ell}^+,z_{\ell}^+)$ be a solution of subproblem \reff{qp1a1}--\reff{qp1a5} at $(x_{\ell},y_{\ell},z_{\ell})$ with parameter $\rho_k$ and $B_{\rho}=B_{\ell}$,
then $y_{\ell}^+=|h(x_{\ell})+\na h(x_{\ell})^Td_{\ell}|,~ z_{\ell}^+=\max\{0,-(g(x_{\ell})+\na g(x_{\ell})^Td_{\ell})\}$, and \begin{eqnarray}
\rho_k\na f(x_{\ell})^Td_{\ell}+\hf d_{\ell}^TB_{\ell}d_{\ell}+(\|y_{\ell}^+\|_1\|+\|z_{\ell}^+\|_1)-(\|y_{\ell}\|_1\|+\|z_{\ell}\|_1)\le 0.
\nonumber\end{eqnarray}

 Now we are ready to present our exact penalty SQP algorithm with inner and outer iterations.
\begin{algorithm}\caption{An exact penalty SQP method with inner and outer iterations}\label{alg1}
{\small \alglist
\item[Given] $x_0\in\Re^{n}$, $B_0\in{\cal S}_{++}^{n}$, $\rho_0>0$, $\sigma\in (0,1)$, $\tau\in (0,1)$, and the tolerance $\epsilon$.
    Let $y_0=|h(x_0)|$, $z_0=\max\{0,-g(x_0)\}$. Set $k:=0$.

\item[Step] 1  For given $\rho_k$, solve the smooth feasible optimization \reff{prob2-1}--\reff{prob2-5} by SQP method approximately.

{\bf Step 1.0} Set $\ell:=0$, $(x_{\ell},y_{\ell},z_{\ell})=(x_k,y_k,z_k)$, $B_{\ell}=B_k$.

{\bf Step 1.1} Solve the QP subproblem \reff{qp1a1}--\reff{qp1a5} to derive $(d_{\ell},y_{\ell}^+,z_{\ell}^+)$ and Lagrange multipliers $(u_{\ell+1},v_{\ell+1},s_{\ell+1},t_{\ell+1})$, where $(x,y,z)=(x_{\ell},y_{\ell},z_{\ell})$, $B=B_{\ell}$, $\rho=\rho_k$.

{\bf Step 1.2} Evaluate $r_{\ell}:=(\|y_{\ell}^+\|_1+\|z_{\ell}^+\|_1)-(\|y_{\ell}\|_1+\|z_{\ell}\|_1)$. If $r_{\ell}\ge-\epsilon$, set $(x_{k+1},y_{k+1},z_{k+1})=(x_{\ell},y_{\ell},z_{\ell})$, $\mu_{k+1}=(v_{\ell+1}-u_{\ell+1})$, $\lambda_{k+1}=s_{\ell+1}$, $B_{k+1}=B_{\ell}$, stop the inner SQP algorithm and go to Step 2.

{\bf Step 1.3} Compute $P_{\rho_k}^{'}(x_{\ell};d_{\ell}):=\rho_{k}\na f(x_{\ell})^Td_{\ell}+r_{\ell}$. Select a nonnegative integer $\imath$ as small as possible such that \begin{eqnarray}
P_{\rho_{k}}(x_{\ell}+\tau^{\imath}d_{\ell})-P_{\rho_{k}}(x_{\ell})\le\sigma\tau^{\imath}P_{\rho_{k}}^{'}(x_{\ell};d_{\ell}) \label{stepsize}\end{eqnarray}
and let $\alp_{\ell}=\tau^{\imath}$.

{\bf Step 1.4} Let $x_{\ell+1}=x_{\ell}+\alp_{\ell}d_{\ell},~y_{\ell+1}=|h(x_{\ell+1})|,~z_{\ell+1}=\max\{0,-g(x_{\ell+1})\}$.
Update $B_{\ell}$ to $B_{\ell+1}$, set $\ell:=\ell+1$ and go to Step 1.1.

\item[Step] 2 ({\bf Update the parameter}). Let $d_{k+1}$ be a solution of the regularized linear $\ell_1$ minimization \begin{eqnarray}\min_d~\hf d^TB_{k+1}d+c'(x_{k+1};d), \label{LS}\end{eqnarray}
$y_{k+1}^+=|h(x_{k+1})+\na h(x_{k+1})^Td_{k+1}|$, $z_{k+1}^+=\max\{0,-(g(x_{k+1})+\na g(x_{k+1})^Td_{k+1})\}$, and $r_{k+1}=(\|y_{k+1}^+\|_1+\|z_{k+1}^+\|_1)-(\|y_{k+1}\|_1+\|z_{k+1}\|_1)$. If $\|d_{k+1}\|_{\infty}>\epsilon$, select a nonnegative integer $\imath$ as small as possible such that $x_{k+1}^0=x_{k+1}+\alp_{k+1}d_{k+1}$ with $\alp_{k+1}=\tau^{\imath}$ satisfies \begin{eqnarray}
c(x_{k+1}^0)-c(x_{k+1})\le\sigma\tau^{\imath}r_{k+1}, \label{stepsize1}\end{eqnarray}
and update $\rho_k$ to $\rho_{k+1}=
\min\{0.01\rho_k,\frac{c(x_{k+1})-c(x_{k+1}^0)}{f(x_{k+1}^0)-f(x_{k+1})}\}$ provided
$\rho_kf(x_{k+1}^0)+c(x_{k+1}^0)>\rho_kf(x_{k+1})+c(x_{k+1})$, otherwise $\rho_{k+1}=
\min\{0.1\rho_k,\rho_k^{1.5}\}$. Update $B_{k+1}$, set $(x_{k+1},y_{k+1},z_{k+1})=(x_{k+1}^0,y_{k+1}^0,z_{k+1}^0)$. Else if $\|d_{\ell}\|_{\infty}>\epsilon$, set $\rho_{k+1}=\min\{0.01\rho_k,\rho_k^{1.5}\}$. 

\item[Step 3] If $\max\{\|d_{k+1}\|_{\infty},\|d_{\ell}\|_{\infty}\}\le\epsilon$, stop the algorithm. Otherwise, set $k:=k+1$ and go to Step 1.

\eli}
\end{algorithm}

Our algorithm consists of inner iterations (Step 1) and Outer iterations (Step 2), where the former approximately solves the smooth equivalent of the exact penalty problem \reff{prob1L1}, and the latter updates the penalty parameter and provides suitable initial points for inner iterations based on feedback of the former. The whole algorithm will terminate at Step 3 with $\|d_{\ell}\|_{\infty}$ and $\|d_k\|_{\infty}$ sufficiently small.
The inner iterations terminate when $r_{\ell}\ge-\epsilon$. This condition will be satisfied as $\|d_{\ell}\|$ is small enough, while it admits the termination of inner algorithm when $\|d_{\ell}\|$ is not very small but the constraint violations could not be reduced apparently.

Our algorithm has similarity to the existing SQP methods with rapid detection of infeasibility \cite{BurCuW14,ByrCuN10} in that it is proposed based an exact penalty function and uses line search methods along directions obtained by solving QP subproblems. However, they are distinct in the way for updating the penalty parameter, which brings about the distinct global and local convergence results.

\section{Global convergence}

Our algorithm consists of the inner loop and the outer loop, in which the inner loop is an SQP algorithm for solving a smooth feasible optimization \reff{prob2-1}--\reff{prob2-5} with a fixed penalty parameter and the outer algorithm updates the penalty parameter based on the feedback of inner algorithm. In particular, since MFCQ always holds, every local minimizer of optimization \reff{prob2-1}--\reff{prob2-5} is a KKT point and the Lagrange multipliers are bounded by the parameter.

\subsection{Global convergence of the inner SQP}
In this subsection, we analyze convergence of the inner SQP algorithm.
It is assumed that the inner SQP algorithm does not terminate in finite number of iterations. Our analysis will show that $d_{\ell}\to 0$ as $\ell\to\infty$ under this supposition, which further implies $r_{\ell}\to 0$ and the finite termination of inner loop with $r_{\ell}\ge-\epsilon$. Based on the supposition, the infinite sequences $\{x_{\ell}\}$, $\{y_{\ell}\}$, $\{z_{\ell}\}$, $\{B_{\ell}\}$, $\{d_{\ell}\}$ are generated. The fixed penalty parameter is $\rho_k$.
Before doing our global analysis, we need the following blanket assumptions.

\begin{ass}\label{ass1} (1) The iterative sequence $\{x_{\ell}\}$ is in a open bounded set of $\Re^n$ (which may depend on $\rho_k$).
(2) The sequence of approximate matrices $\{B_{\ell}\}$ is uniformly bounded and symmetric positive definite for all $\ell$.
\end{ass}

The conditions (1) and (2) in \refa{ass1} are quite general conditions extensively used in convergence analysis for nonlinear optimization.
\refa{ass1} (2) is a moderate and easy satisfied condition since, as an example, the approximate Hessian $B_{\ell}$ can be simply chosen to be a multiple of the identity matrix.

The first two results show that every QP subproblem for the exact penalty optimization is well-behaved, and the line search procedure is well-defined.
\begin{lemma} Under \refa{ass1}, for any given $\rho_k$, the solution $(d_{\ell},y_{\ell}^+,z_{\ell}^+)$ of the QP subproblem is unique and the associated vectors of Lagrange multipliers $(u_{\ell+1},v_{\ell+1},s_{\ell+1},t_{\ell+1})$ are bounded in $[0,1]$. Furthermore, if all gradients $\na h_i(x_{\ell}), ~i\in\{i|h_i(x_{\ell})=0\}$ and $\na g_i(x_{\ell}), ~i\in\{i| g_i(x_{\ell})=0\}$ are linearly independent, then the vector of multipliers $(u_{\ell+1},v_{\ell+1},s_{\ell+1},t_{\ell+1})$ is also unqiue. \end{lemma}\begin{proof}
For simplicity of statement, let \begin{eqnarray}
M_{\ell}(d)=\dd\dd\rho_k\na f(x_{\ell})^Td+\hf d^TB_{\ell}d+(\|h(x_{\ell})+\na h(x_{\ell})^Td\|_1 \nonumber\\
\dd\dd+\|\max\{0,-(g(x_{\ell})+\na g(x_{\ell})^Td)\}\|_1). \nonumber\end{eqnarray}
We first prove that $M_{\ell}(d)$ is a strongly convex function. By some calculation, for any $d_1\in\Re^n$, $d_2\in\Re^n$, and any $\alpha\in[0,1]$, one has \begin{eqnarray}
M_{\ell}(\alpha d_1+(1-\alpha)d_2)\le\alpha M_{\ell}(d_1)+(1-\alpha)M_{\ell}(d_2)-\hf\alpha(1-\alpha)(d_1-d_2)^TB_{\ell}(d_1-d_2) \nonumber\end{eqnarray} since
$\|h(x_{\ell})+\na h(x_{\ell})^T(\alpha d_1+(1-\alpha)d_2)\|_1
\le\alpha\|h(x_{\ell})+\na h(x_{\ell})^Td_1\|_1+(1-\alpha)\|h(x_{\ell})+\na h(x_{\ell})^Td_2\|_1$ and
\begin{eqnarray}
\dd\dd \|\max\{0,-(g(x_{\ell})+\na g(x_{\ell})^T(\alpha d_1+(1-\alpha)d_2))\}\|_1 \nonumber\\
\dd\dd \le\alpha\|\max\{0,-(g(x_{\ell})+\na g(x_{\ell})^Td_1)\}\|_1+(1-\alpha)\|\max\{0,-(g(x_{\ell})+\na g(x_{\ell})^Td_2)\}\|_1.\nonumber\end{eqnarray}
Then the strong convexity of $M_{\ell}(d)$ follows from the uniform symmetric positive definiteness of $B_{\ell}$ (see \refa{ass1} (2)).

The fact that $M_{\ell}(d)$ is strongly convex implies that $d_{\ell}$ is unique, which further implies that \begin{eqnarray}
y^+_{\ell}=|h(x_{\ell})+\na h(x_{\ell})^Td_{\ell}|\quad\hbox{and}\quad z^+_{\ell}=|\max\{0,-(g(x_{\ell})+\na g(x_{\ell})^Td_{\ell})\}| \nonumber\end{eqnarray}
are unique. In addition, it follows from the KKT conditions of QP subproblems that \begin{eqnarray}
u_{\ell+1,i}+v_{\ell+1,i}=1,~\forall i=1,\ldots,m_{\cal E};~s_{\ell+1,i}+t_{\ell+1,i}=1,~\forall i=1,\ldots,m_{\cal I}, \nonumber\end{eqnarray}
and $(u_{\ell+1},v_{\ell+1},s_{\ell+1},t_{\ell+1})\ge 0$. Thus, $\|(u_{\ell+1},v_{\ell+1},s_{\ell+1},t_{\ell+1})\|_{\infty}\le1$.

If the vectors of Lagrange multipliers $(u_{\ell+1},v_{\ell+1},s_{\ell+1},t_{\ell+1})$ are not unique, and an other vectors $(u_{\ell+1}^{'},v_{\ell+1}^{'},s_{\ell+1}^{'},t_{\ell+1}^{'})\ne(u_{\ell+1},v_{\ell+1},s_{\ell+1},t_{\ell+1})$ also satisfy the KKT conditions of QP subproblems, then \begin{eqnarray}
\na h(x_{\ell})(u_{\ell+1}^{'}-v_{\ell+1}^{'}-u_{\ell+1}+v_{\ell+1})-\na g(x_{\ell})(s_{\ell+1}^{'}-s_{\ell+1})=0. \nonumber \end{eqnarray}
Noting that $u_{\ell+1,i}^{'}-v_{\ell+1,i}^{'}-u_{\ell+1,i}+v_{\ell+1,i}=0$ for $i\notin\{i|h_i(x_{\ell})=0\}$ and $s_{\ell+1,i}^{'}-s_{\ell+1,i}=0$ for $i\notin\{i| g_i(x_{\ell})=0\}$,
the preceding equation contradicts that the gradients $\na h_i(x_{\ell}), ~i\in\{i|h_i(x_{\ell})=0\}$ and $\na g_i(x_{\ell}), ~i\in\{i| g_i(x_{\ell})=0\}$ are linearly independent, which completes our proof. \end{proof}

\begin{lemma}\label{alp}  Under \refa{ass1}, there exists a scalar $\beta_k\in (0,1]$ independent of $\ell$ such that for $\alpha_{\ell}\in (0,\beta_k]$, \reff{stepsize} always holds. Thus, $\{\alpha_{\ell}\}$ is bounded away from zero.
\end{lemma}\begin{proof} We firstly prove that $\{d_{\ell}\}$ is bounded by contradiction. Due to \begin{eqnarray}
\rho_k\na f(x_{\ell})+B_{\ell}d_{\ell}+\na h(x_{\ell})(u_{\ell+1}-v_{\ell+1})-\na g(x_{\ell})s_{\ell+1}=0, \nonumber\end{eqnarray}
and the boundedness of $B_{\ell}$, $\na f(x_{\ell})$, $\na h(x_{\ell})$, and $(u_{\ell+1},v_{\ell+1},s_{\ell+1},t_{\ell+1})$, if $\{d_{\ell}\}$ is unbounded, dividing $\|d_{\ell}\|$ and taking the limit $\ell\to\infty$ on both sides of the preceding equation will result in $\lim_{\ell\to\infty}B_{\ell}\frac{d_{\ell}}{\|d_{\ell}\|}=0$, which contradicts the uniform positive definiteness of $B_{\ell}$.

For every $\ell\ge 0$, one has $P^{'}_{\rho_k}(x_{\ell};d_{\ell})\le-\hf d_{\ell}^TB_{\ell}d_{\ell}\le-\hf\lambda_{\min}(B_{\ell})\|d_{\ell}\|^2$ and \begin{eqnarray} \dd\dd P_{\rho_k}(x_{\ell}+\alp d_{\ell})-P_{\rho_k}(x_{\ell})-\sigma\alp P^{'}_{\rho_k}(x_{\ell};d_{\ell}) \nonumber\\
\dd\dd=\rho_k f(x_{\ell}+\alp d_{\ell})-\rho_k f(x_{\ell})-\sigma\alp\rho_k\na f(x_{\ell})^Td_{\ell} \nonumber\\
\dd\dd\quad+(\|h(x_{\ell}+\alp d_{\ell})\|_1-\|h(x_{\ell})\|_1-\sigma\alp(\|h(x_{\ell})+\na h(x_{\ell})^Td_{\ell}\|_1-\|h(x_{\ell})\|_1)) \nonumber\\
\dd\dd\quad+(\|\max\{0,-g(x_{\ell}+\alp d_{\ell})\}\|_1-\|\max\{0,-g(x_{\ell})\}\|_1 \nonumber\\
\dd\dd\quad\quad\quad-\sigma\alp(\|\max\{0,-(g(x_{\ell})+\na g(x_{\ell})^Td_{\ell})\}\|_1-\|\max\{0,-g(x_{\ell})\}\|_1)) \nonumber\\
\dd\dd\le(1-\sigma)\alp P^{'}_{\rho_k}(x_{\ell};d_{\ell})+O(\alp^2\|d_{\ell}\|^2). \label{230405a}\end{eqnarray}
Then, due to the boundedness of $\|d_{\ell}\|$ and $\sigma\in (0,1)$, there exists a scalar $\beta_k\in (0,1]$ independent of inner iterations such that the left side of \reff{230405a} is non-positive for all $\alpha\in (0,\beta_k]$. Consequently, the result follows immediately.
\end{proof}

The next result indicates the finite termination of inner iterations.
\begin{theorem} Under \refa{ass1}, $\{P_{\rho_k}(x_{\ell})\}$ is lower bounded. Then one has $\lim_{\ell\to\infty} \|d_{\ell}\|=0$. That is, for any given $\epsilon>0$, $r_{\ell}\ge-\epsilon$ will be satisfied in a finite number of iterations. \end{theorem}\begin{proof}
Note that the sequence $\{P_{\rho_k}(x_{\ell})\}$ is monotonically non-increasing. Thus, one has $\lim_{\ell\to\infty}P_{\rho_k}(x_{\ell})=P_k^*$ with $P_k^*$ being finite or $-\infty$. Since $f$ is lower bounded, $P_{\rho_k}(x_{\ell})$ is lower bounded, which rules the case $P_k^*=-\infty$ out. By taking the limit $\ell\to\infty$ on both hand-sides of \reff{stepsize} and noting the result of \refl{alp} and the inequality $P^{'}_{\rho_k}(x_{\ell};d_{\ell})\le-\hf\lambda_{\min}(B_{\ell})\|d_{\ell}\|^2$, one has the limit $\lim_{\ell\to\infty} \|d_{\ell}\|=0$ which further implies $\lim_{\ell\to\infty} |r_{\ell}|=0$. Our proof is completed. \end{proof}

\subsection{Global convergence of the whole algorithm}
In this subsection, we prove that for small $\rho_k$, $x_k$ can be either an approximate KKT point or an approximate DL-stationary point of the original problem. Otherwise, $\rho_k\to 0$ and there exists a cluster point of $\{x_k\}$ which is either a singular stationary point or a DZ-stationary point of problem \reff{prob1-1}.

Firstly, by our definitions, it is easy to know that if $\|y_k\|_1=0$ and $\|z_k\|_1=0$, then $x_k$ is a KKT point of the original optimization problem \reff{prob1-1}. When $\|y_k\|_1+\|z_k\|_1\ne 0$, $x_k$ may be a DL-stationary point of the original optimization problem \reff{prob1-1} when $d_k=0$.
\begin{lemma}\label{s4.2.2} If $x_k$ is a KKT point of the smooth feasible problem \reff{prob2-1}--\reff{prob2-5} and $y_k=z_k=0$, then $x_k$ is a KKT point of the nonlinear optimization problem \reff{prob1-1}. \end{lemma}\begin{proof}
Due to $y_k=z_k=0$, $x_k$ is feasible to the problem \reff{prob1-1}. Then the KKT conditions \reff{kkt2.9-1}--\reff{kkt2.9-8}
of problem \reff{prob2-1}--\reff{prob2-5} are reduced to those for problem \reff{prob1-1}. \end{proof}

\begin{lemma} If for some $k\ge 0$, the inner iterations $\{x_{\ell}\}$ of the SQP algorithm converge to $x_{k+1}$ which is infeasible to the problem \reff{prob1-1}, and $d_{k+1}=0$ is a solution of the regularized linear $\ell_1$ minimization subproblem \reff{LS}, then $x_{k+1}$ is a DL-stationary point of the original optimization problem \reff{prob1-1}.  \end{lemma}\begin{proof}
The fact $d_{k+1}=0$ shows $x_{k+1}$ is a stationary point of the $\ell_1$ minimization problem of constraint violations. Thus by \refd{def1}, $x_{k+1}$ is a DL-stationary point of the original optimization problem \reff{prob1-1}. \end{proof}

We are left to show the results on $\rho_k$ to be sufficiently small, which will be given in the next lemma.
\begin{lemma} Let $\{x_k\}$ be an infinite sequence generated by the outer algorithm. Then $k\to\infty$ and $\rho_k\to 0$. Suppose that ${\cal K}$ is an infinite set of indices and $\{x_k| k\in{\cal K}\}$ is any convergent subsequence of $\{x_k\}$, $\|y_k\|_1+\|z_k\|_1\to 0$ as $k\in{\cal K}$ and $k\to\infty$, then the limit point of $\{x_k| k\in{\cal K}\}$ is a singular stationary point of problem \reff{prob1-1}. Otherwise, $\|y_k\|_1+\|z_k\|_1\not\to 0$ as $k\in{\cal K}$ and $k\to\infty$ and every cluster point of $\{x_k| k\in{\cal K}\}$ is a DZ-stationary point of problem \reff{prob1-1}.
\end{lemma}\begin{proof} For $k\in{\cal K}$, the limit $\|y_k\|_1+\|z_k\|_1\to 0$ as $k\to\infty$ implies that the limit point of $\{x_k| k\in{\cal K}\}$ is a feasible point of problem \reff{prob1-1}. Combining with the KKT conditions of problem \reff{prob2-1}--\reff{prob2-5} and $\rho_k\to 0$ results in that the stationary conditions given in \refd{def0} (2) hold at the limit point. Thus, the limit point of $\{x_k| k\in{\cal K}\}$ is a singular stationary point of problem \reff{prob1-1}.

Similarly, if $\|y_k\|_1+\|z_k\|_1\not\to 0$ as $k\in{\cal K}$ and $k\to\infty$, then any cluster point of $\{x_k| k\in{\cal K}\}$ is infeasible to the problem \reff{prob1-1}, which brings about DZ-stationary points of problem \reff{prob1-1} since both the KKT conditions of problem \reff{prob2-1}--\reff{prob2-5} and the limit $\rho_k\to 0$ imply the conditions given in \refd{def1} (2). \end{proof}

In summary, we have the main global convergence results on \refal{alg1} as follows.
\begin{theorem} Suppose that $\{(x_k,y_k,z_k)\}$ is a sequence generated by the outer algorithm. If $(x_k,y_k,z_k)$ is the terminating point with $\rho_k$ far from zero, then $x_k$ is an approximate KKT point of problem \reff{prob1-1} provided both $\|y_k\|$ and $\|z_k\|$ small enough, otherwise is an approximate DL-stationary point of problem \reff{prob1-1}. If $(x_k,y_k,z_k)$ is the terminating point with $\rho_k$ close to zero enough, then $x_k$ is an approximate singular stationary point of problem \reff{prob1-1} provided both $\|y_k\|$ and $\|z_k\|$ small enough, otherwise is an approximate DZ-stationary point of problem \reff{prob1-1}. \end{theorem}\begin{proof}
The results follow from the preceding three lemmas of this subsection immediately.  \end{proof}

\section{Local convergence}

We have proved in the preceding section that under very mild conditions, our algorithm may either terminate finitely at certain outer iteration with $\rho_k$ far from zero or generate a sequence of outer iterations with $\rho_k\to 0$. When the outer iteration terminates with $\rho_k$ far from zero, the inner iterations may converge to either a KKT point or a DL-stationary point; otherwise, the outer iterations may converge to either a singular stationary point or a DZ-stationary point.
In Section 6.1, it is proved that under suitable assumptions, our algorithm can be quadratically convergent to a KKT point when the original problem is feasible, and can be quadratically convergent to a DL-stationary point when it is infeasible (that is, {\sl our algorithm may rapidly detect infeasibility of the problem without requiring $\rho_k$ close to zero}). In Section 6.2, our algorithm can be rapidly convergent to a DZ-stationary point as $\rho_k\to 0$, which is similar to the results of \cite{BurCuW14,ByrCuN10}. Like all local convergence analysis, we need to suppose that all functions $f(x)$, $g(x)$ and $h(x)$ are twice differentiable, and their second derivatives are Lipschitz continuous around the limit $x^*$.

\subsection{Rapid convergence to either a KKT point or a DL-stationary point}
This subsection focuses on rapid convergence of the inner algorithm for fixed $\rho_k>0$.
It is well known that SQP can be locally quadratically/superlinearly convergent when the sequence of Hessians of subproblems is suitably approximated to the exact Lagrangian Hessian, LICQ and the strict complementarity conditions hold at the solution. Thus, we will examine these conditions on problem \reff{prob2-1}--\reff{prob2-5}.

\begin{ass}\label{ass2a} Let $\{(x_{\ell},y_{\ell},z_{\ell})\}$ be the iterative sequence generated by the inner algorithm, and $\{(u_{\ell},v_{\ell},s_{\ell},t_{\ell})\}$ is the associated sequence of Lagrangian multipliers. There hold $(x_{\ell},y_{\ell},z_{\ell})\to (x^*,y^*,z^*)$ and $(u_{\ell},v_{\ell},s_{\ell},t_{\ell})\to(u^*,v^*,s^*,t^*)$ as $\ell\to\infty$, where $(x^*,y^*,z^*)$ is a KKT point of problem \reff{prob2-1}--\reff{prob2-5} in which $x^*$ is either a KKT point (correspondingly, $y^*+z^*=0$) or a DL-stationary point ($y^*+z^*>0$) of the original problem \reff{prob1-1}, $(u^*,v^*,s^*,t^*)$ is the associated Lagrange multiplier vector.
\end{ass}

Under \refa{ass2a}, all components of the Lagrange multiplier vectors $u^*$, $v^*$, $s^*$, $t^*$ are in the interval $[0,1]$.
 Let ${\cal I}_h^{\ell}=\{i\in{\cal E}| h_i(x_{\ell})=0\}$, ${\cal I}_g^{\ell}=\{i\in{\cal I}| g_i(x_{\ell})=0\}$, ${\cal I}_h^{*}=\{i\in{\cal E}| h_i(x^*)=0\}$, ${\cal I}_g^{*}=\{i\in{\cal I}| g_i(x^*)=0\}$, and $\mu_{\ell}=\frac{v_{\ell}-u_{\ell}}{\rho_k}$, $\lambda_{\ell}=\frac{s_{\ell}}{\rho_k}$, $\mu^*=\frac{v^*-u^*}{\rho_k}$, $\lambda^*=\frac{s^*}{\rho_k}$. Then $\mu_{\ell}\to\mu^*$ and $\lambda_{\ell}\to\lambda^*$ as $\ell\to\infty$.

\begin{ass}\label{ass2c1} \ \\
(1) The gradients $\na h_i(x^*), i\in{\cal I}_h^*$ and $\na g_i(x^*), i\in{\cal I}_g^*$ are linearly independent;\\
(2) The multipliers $u_i^*\in (0,1)$ and $v_i^*\in (0,1)$, $\forall~i\in{\cal I}_h^*$, and $s_i^*\in (0,1)$, $\forall~i\in{\cal I}_g^*$; \\
(3) The Lagrangian Hessian satisfies $$d^T\na^2_{xx}L(x^*,\mu^*,\lambda^*)d\ge\gamma\|d\|^2,$$ $\forall d\in\{d| \na h_i(x^*)^Td=0,~i\in{\cal I}_h^*; \na g_i(x^*)^Td=0,~i\in{\cal I}_g^*\}$, where $L(x,\mu,\lambda)=f(x)-\mu^Th(x)-\lambda^Tg(x)$, and $\gamma>0$ is a constant.
\end{ass}

\refa{ass2c1} (1) holds if and only if LICQ holds with problem \reff{prob2-1}--\reff{prob2-5}. Since $u_i^*+v_i^*=1$ for all $i\in{\cal E}$, \refa{ass2c1} (2) implies the strict complementarity conditions for problem \reff{prob2-1}--\reff{prob2-5}. Correspondingly, $|\rho_k\mu^*_i|<1,~i\in{\cal I}_h^*$ and $\rho_k\lambda^*_i\in(0,1),~i\in{\cal I}_g^*$ for the original problem.
The following results show that under Assumptions \ref{ass2a} and \ref{ass2c1}, $d_{\ell}$
can be a superlinearly or quadratically convergent step.

\begin{theorem}\label{thms1} Suppose that Assumptions \ref{ass2a} and \ref{ass2c1} hold.

(1) If $\|(B_{\ell}-\na^2_{xx} L(x^*,\mu^*,\lambda^*))d_{\ell}\|=o(\|d_{\ell}\|)$, then
\begin{eqnarray}\lim_{\ell\to\infty}\frac{\|x_{\ell}+d_{\ell}-x^*\|}{\|x_{\ell}-x^*\|}=0.
\label{local1}\end{eqnarray}

(2) If $B_{\ell}=\na^2_{xx}L(x_{\ell},\mu_{\ell},\lambda_{\ell})$, then $\|x_{\ell}+d_{\ell}-x^*\|=O(\|x_{\ell}-x^*\|^2).$
\end{theorem}\begin{proof} Let $A_*=[\na h_{{\cal I}_h^*}(x^*)~\na g_{{\cal I}_g^*}(x^*)]$, $A_{\ell}=[\na h_{{\cal I}_h^{\ell}}(x_{\ell})~\na g_{{\cal I}_g^{\ell}}(x_{\ell})]$, $P_*=I-A_*(A_*^TA_*)^{-1}A_*^T$ and $P_{\ell}=I-A_{\ell}(A_{\ell}^TA_{\ell})^{-1}A_{\ell}^T$, where $I$ is the
$n\times n$ identity matrix. Consider the system \begin{eqnarray} \dd\dd
\left[\begin{array}{c} P_*\na_{xx}^2L(x^*,\mu^*,\lambda^*) \\ A_*^T\end{array}\right]d=0. \label{080915c}\end{eqnarray} Let $d^*\in\Re^n$ be any one
of its solutions. If $d^*\ne 0$, then \begin{eqnarray}
(d^*)^TP_*\na_{xx}^2L(x^*,\mu^*,\lambda^*)d^*=0, \quad A_*^Td^*=0, \nonumber\end{eqnarray} so we have $
(d^*)^T\na_{xx}^2L(x^*,\mu^*,\lambda^*)d^*=0, $ which contradicts
\refa{ass2c1} (3). This contradiction shows that the coefficient
matrix of the system \reff{080915c} has full column rank.
Therefore, by Assumptions \ref{ass2a} and \ref{ass2c1}, for all sufficiently large ${\ell}$,
the matrix \begin{eqnarray} \left[\begin{array}{c} P_{\ell}\na_{xx}^2L(x^*,\mu^*,\lambda^*) \\
A_{\ell}^T\end{array}\right] \nonumber\end{eqnarray} is of full column rank.

Since $\na f(x_{\ell})+B_{\ell}d_{\ell}-\na h_{{\cal I}_h^{\ell}}(x_{\ell})\mu_{{\cal I}_h^{\ell}}-\na g_{{\cal I}_g^{\ell}}(x_{\ell})\lambda_{{\cal I}_g^{\ell}}=0,$ there holds \begin{eqnarray}
 P_{\ell}B_{\ell}d_{\ell}
\dd\dd =-P_{\ell}(\na f(x_{\ell})-\na h_{{\cal I}_h^{\ell}}(x_{\ell})\mu_{{\cal I}_h^{\ell}}-\na g_{{\cal I}_g^{\ell}}(x_{\ell})\lambda_{{\cal I}_g^{\ell}}) \nonumber\\
\dd\dd =-P_{\ell}(\na f(x_{\ell})-\na h_{{\cal I}_h^{\ell}}(x_{\ell})\mu_{{\cal I}_h^*}-\na g_{{\cal I}_g^{\ell}}(x_{\ell})\lambda_{{\cal I}_g^*}) \nonumber\\
\dd\dd
=-P_{\ell}\na_{xx}^2L(x^*,\mu^*,\lambda^*)(x_{\ell}-x^*)+{O}(\|x_{\ell}-x^*\|^2).
\nonumber \end{eqnarray} Thus, \begin{eqnarray}
\dd\dd P_{\ell}(B_{\ell}-\na_{xx}^2L(x^*,\mu^*,\lambda^*))d_{\ell} \nonumber\\
\dd\dd =-P_{\ell}\na_{xx}^2L(x^*,\mu^*,\lambda^*)
(x_{\ell}+d_{\ell}-x^*)+{O}(\|x_{\ell}-x^*\|^2).
\label{080915a}\end{eqnarray}
Thanks to $h_{{\cal I}_h^{\ell}}(x^*)=0$ and $h_{{\cal I}_h^{\ell}}(x_{\ell})=h_{{\cal I}_h^{\ell}}(x_{\ell})-h_{{\cal
                        I}_h^{\ell}}(x^*) =\na h_{{\cal
                        I}_h^{\ell}}(x^*)^T(x_{\ell}-x^*)+{O}(\|x_{\ell}-x^*\|^2)$, and $g_{{\cal I}_g^{\ell}}(x^*)=0$ and $g_{{\cal I}_g^{\ell}}(x_{\ell})=g_{{\cal I}_g^{\ell}}(x_{\ell})-g_{{\cal
                        I}_g^{\ell}}(x^*) =\na g_{{\cal
                        I}_g^{\ell}}(x^*)^T(x_{\ell}-x^*)+{O}(\|x_{\ell}-x^*\|^2)$, there holds
\begin{eqnarray} A_{{\ell}}^T(x_{\ell}+d_{\ell}-x^*)=\left[\begin{array}{c}
h_{{\cal I}_h^{\ell}}(x_{\ell}) \\
g_{{\cal I}_g^{\ell}}(x^*)\end{array}\right]+A_{{\ell}}^Td_{\ell}+{O}(\|x_{\ell}-x^*\|^2). \label{080915b}\end{eqnarray}

Putting \reff{080915a} and \reff{080915b} together in a matrix, we can
obtain \begin{eqnarray}  \dd\dd\left[\begin{array}{c} P_{\ell}\na_{xx}^2L(x^*,\mu^*,\lambda^*) \\
A_{{\ell}}^T\end{array}\right](x_{\ell}+d_{\ell}-x^*)\nonumber\\
\dd\dd =\left[\begin{array}{c} -P_{\ell}(B_{\ell}-\na_{xx}^2L(x^*,\mu^*,\lambda^*))d_{\ell} \\
h_{{\cal I}_h^{\ell}}(x_{\ell})+\na h_{{\cal I}_h^{\ell}}(x_{\ell})^Td_{\ell} \\
g_{{\cal I}_g^{\ell}}(x_{\ell})+\na g_{{\cal I}_g^{\ell}}(x_{\ell})^Td_{\ell}
\end{array}\right]+{O}(\|x_{\ell}-x^*\|^2), \label{080915d}\end{eqnarray}
where the coefficient matrix has previously been proved to be of
full column rank for all sufficiently large ${\ell}$.

(1) If $\|(B_{\ell}-\na^2_{xx} L(x^*,\mu^*,\lambda^*))d_{\ell}\|=o(\|d_{\ell}\|)$, then
\begin{eqnarray}
\left\|\left[\begin{array}{c} -P_{\ell}(B_{\ell}-\na_{xx}^2L(x^*,\mu^*,\lambda^*))d_{\ell} \\
h_{{\cal I}_h^{\ell}}(x_{\ell})+\na h_{{\cal I}_h^{\ell}}(x_{\ell})^Td_{\ell} \\
g_{{\cal I}_g^{\ell}}(x_{\ell})+\na g_{{\cal I}_g^{\ell}}(x_{\ell})^Td_{\ell}
\end{array}\right]\right\|={o}(\|d_{\ell}\|). \nonumber\end{eqnarray} Hence, by \reff{080915d} and since
$x_{\ell}\to x^*$ as ${\ell}\to\infty$, there holds \begin{eqnarray}
\lim_{\ell\to\infty}\frac{\|x_{\ell}+d_{\ell}-x^*\|}{\|x_{\ell}-x^*\|}=\lim_{{\ell}\to\infty}\frac{{o}(\|d_{\ell}\|)}{\|x_{\ell}-x^*\|},
\nonumber \end{eqnarray} which shows  \begin{eqnarray}
\lim_{{\ell}\to\infty}\frac{\|d_{\ell}\|}{\|x_{\ell}-x^*\|}=1 \label{081104a}\end{eqnarray} and
thus \reff{local1} follows immediately.

(2) If $B_{\ell}=\na^2_{xx}L(x_{\ell},\mu_{\ell},\lambda_{\ell})$,
then \begin{eqnarray}
\left\|\left[\begin{array}{c} -P_{\ell}(B_{\ell}-\na_{xx}^2L(x^*,\mu^*,\lambda^*))d_{\ell} \\
h_{{\cal I}_h^{\ell}}(x_{\ell})+\na h_{{\cal I}_h^{\ell}}(x_{\ell})^Td_{\ell} \\
g_{{\cal I}_g^{\ell}}(x_{\ell})+\na g_{{\cal I}_g^{\ell}}(x_{\ell})^Td_{\ell}
\end{array}\right]\right\|=O(\|x_{\ell}-x^*\|\|d_{\ell}\|), \nonumber\end{eqnarray}
which implies that \reff{081104a} holds. Hence, it follows from \reff{080915d} that \begin{eqnarray}
\lim_{\ell\to\infty}\frac{\|x_{\ell}+d_{\ell}-x^*\|}{\|x_{\ell}-x^*\|^2}=\lim_{{\ell}\to\infty}\frac{{O}(\|d_{\ell}\|)}{\|x_{\ell}-x^*\|}+O(1)=O(1),
\nonumber\end{eqnarray}
which completes our proof. \end{proof}

\subsection{Rapid convergence to a DZ-stationary point}
We have known that convergence to a DZ-stationary point can happen only when $k\to\infty$ and $\rho_k\to 0$.
For convenience of statement, we denote $w^*=(x^*,\mu^*,\la^*)$, $\hat w^*=(x^*,u^*,v^*,s^*,t^*)$, and $w_k=(x_k,\mu_k,\la_k)$, $\hat w_k=(x_k,u_k,v_k,s_k,t_k)$ for all $k\ge 0$. When applying to optimization with only inequality constraints, the approach for local analysis in this subsection is almost the same as that of \cite{ByrCuN10}. Thus, it can be thought of as an extension of the analysis of \cite{ByrCuN10} to general constrained optimization.
The following assumptions are similar to those of \cite{ByrCuN10}.

\begin{ass}\label{ass2b} \ \\
(1) $w_k\rightarrow  w^*$ and $\rho_k\rightarrow  0$ as $k\rightarrow \infty$; \\
(2) The gradients $\na h_i(x^*), i\in{\cal I}_h^*$ and $\na g_i(x^*), i\in{\cal I}_g^*$ are linearly independent;\\
(3) The multipliers $u_i^*\in (0,1)$ and $v_i^*\in (0,1)$, $\forall~i\in{\cal I}_h^*$, and $s_i^*\in (0,1)$, $\forall~i\in{\cal I}_g^*$;  \\
(4) $d^TB^*d>0$ for all $d\ne 0$ such that $\na h_i(x^*)^Td=0$ $\forall i\in{\cal I}_h^*$ and $\na g_i(x^*)^Td=0$ $\forall i\in{\cal I}_g^*$, where $B^*=\sum_{i=1}^{m_{\cal E}}\mu_i^{*}\na^2h_i(x^*)-\sum_{i=1}^{m_{\cal I}}\lambda_i^{*}\na^2g_i(x^*)$, $\mu_i^*=\hbox{sign}(h_i(x^*))$ as $h_i(x^*)\ne 0$ and $\mu_i^*=u_i^*-v_i^*$ for $i\in{\cal I}_h^*$, $\lambda_i^*=1$ as $g_i(x^*)<0$, $\lambda_i^*=0$ as $g_i(x^*)>0$ and $\lambda_i^*=s_i^*$ for $i\in{\cal I}_g^*$.
\end{ass}

\refa{ass2b} (1) and (3) imply that for problems \reff{prob2-1}--\reff{prob2-5}, LICQ holds at $(x_k,y_k,z_k)$ for all sufficiently large $k$. Thus, $(u_k,v_k,s_k,t_k)$ and $(\mu_k,\lambda_k)$ are unique for all sufficiently large $k$. Let
${\cal I}_{h+}^*=\{i| h_i(x^*)>0,~i=1,\ldots,m_{\cal E}\}$, ${\cal I}_{h-}^*=\{i| h_i(x^*)<0,~i=1,\ldots,m_{\cal E}\}$, ${\cal I}_{h}^*=\{i| h_i(x^*)=0,~i=1,\ldots,m_{\cal E}\}$, ${\cal I}_{g-}^*=\{i| g_i(x^*)<0,~i=1,\ldots,m_{\cal I}\}$, ${\cal I}_{g}^*=\{i| g_i(x^*)=0,~i=1,\ldots,m_{\cal I}\}$.
 Consider the system \begin{eqnarray}
\dd\dd\rho_k \na f(x)+\sum_{i\in{\cal I}_{h+}^*}\na h_i(x)-\sum_{i\in{\cal I}_{h-}^*}\na h_i(x)+\sum_{i\in{\cal I}_{h}^*}\mu_i\na h_i(x) \nonumber\\
\dd\dd\quad\quad\quad\quad-\sum_{i\in{\cal I}_{g-}^*}\na g_i(x)-\sum_{i\in{\cal I}_{g}^*}\lambda_i\na g_i(x)=0,  \nonumber\\
\dd\dd h_i(x)=0,~i\in{\cal I}_h^*, \nonumber\\
\dd\dd g_i(x)=0,~i\in{\cal I}_g^*. \nonumber\end{eqnarray}
Let $x_k^{\rho_k}$ be one of its solutions for given value of penalty parameter $\rho_k$. Correspondingly, $w_k^{\rho_k}=(x_k^{\rho_k},\mu_k^{\rho_k},\la_k^{\rho_k})$, $\hat w_k^{\rho_k}=(x_k^{\rho_k},u_k^{\rho_k},v_k^{\rho_k},s_k^{\rho_k},t_k^{\rho_k})$, and like those for $x^*$, we have the index sets ${\cal I}_{h+}^{\rho_k}$, ${\cal I}_{h-}^{\rho_k}$, ${\cal I}_h^{\rho_k}$, ${\cal I}_{g+}^{\rho_k}$, ${\cal I}_{g-}^{\rho_k}$, ${\cal I}_g^{\rho_k}$. The following results can be proved in a similar way to Lemmas 3.3 and 3.5 of \cite{ByrCuN10}.
\begin{lemma}\label{230413a} Under \refa{ass2b}, there exists a $\rho^*>0$ such that the preceding system of equations has a unique solution $w_k^{\rho_k}$ for all $\rho_k\le\rho^*$, and there holds
\begin{equation} \|w_k^{\rho_k}-w^*\|\le {\rho_k}M_1<\epsilon, \label{20140415e1}\end{equation} where $\epsilon>0$ is small enough and
$M_1$ is a positive constant independent of $\rho_k$. Moreover, ${\cal I}_{h+}^{\rho_k}={\cal I}_{h+}^*$, ${\cal I}_{h-}^{\rho_k}={\cal I}_{h-}^*$, ${\cal I}_h^{\rho_k}={\cal I}_h^*$, ${\cal I}_{g+}^{\rho_k}={\cal I}_{g+}^*$, ${\cal I}_{g-}^{\rho_k}={\cal I}_{g-}^*$, ${\cal I}_g^{\rho_k}={\cal I}_g^*$, and \begin{eqnarray}
u_i^{\rho_k}\in (0,1)~\hbox{and}~v_i^{\rho_k}\in (0,1)~\forall~i\in{\cal I}_h^{\rho_k},\quad s_i^{\rho_k}\in (0,1)~\forall~i\in{\cal I}_g^{\rho_k}. \nonumber\end{eqnarray}\end{lemma}
\begin{lemma}\label{230413b} Under \refa{ass2b}, there exist sufficiently small scalars $\rho^*>0$ and $\epsilon>0$ such that for all $\rho_k\le\rho^*$ and $\|w_k-w^*\|<\epsilon$,
\begin{equation} \|w_{k+1}-w_k^{\rho_k}\|\le M_2\|w_k-w_k^{\rho_k}\|, \label{20140415ea}\end{equation}
where we take $x_{k+1}=x_k+d_k$ in $w_{k+1}$, $M_2>0$ is a constant independent of $\rho_k$.
\end{lemma}

In accordance with the preceding two lemmas, one can have the next main result of this subsection.
\begin{theorem} Under \refa{ass2b}, \begin{equation}
\|w_{k+1}-w^*\|=O({\rho_k})+O(\|w_{k}-w^*\|^2). \end{equation}
Therefore, if ${\rho_k}=O(\|w_k-w^*\|)^2$, then the convergence is quadratic; otherwise,
if instead ${\rho_k}=o(\|w_k-w^*\|)$, the convergence is superlinear.
\end{theorem}\begin{proof} The result is straightforward from Lemmas \ref{230413a} and \ref{230413b}. \end{proof}

\section{Numerical experiments}

Our targets in this section are to show that our method is usable and to demonstrate that our theoretical results are achievable.
In this sense, we will not attempt to compare our method with any recognized software, but use our method to solve some small benchmark test examples in the literature for nonlinear programs, for example, \cite{ByrCuN10,DLS17,HocSch81,LDHmm,LiuDHS,LiuSun01}.

We have solved five examples in our numerical experiments, where the first three examples are infeasible and are solved by the methods in Byrd, Curtis and Nocedal \cite{ByrCuN10} and our paper \cite{DLS17} for observing the rapid detection of infeasibility. The fourth and the fifth are feasible examples, where the former is challenging since many effective solvers based on SQP and IPM can only find an infeasible point which has been known in this paper to be a DL-stationary point when starting from some infeasible points, and the latter is one for which the minimizer is a singular stationary point and the linear independence constraint qualification (LICQ) does not hold.

It is noted that, for problem (TP3), our algorithm identifies that the infeasible stationary point is a DL-stationary point, which brings about the rapid detection of infeasibility without driving the penalty parameter $\rho_k$ to zero. This realizes the commentary made by Byrd, Curtis and Nocedal at section 3.2 of their pioneering work \cite{ByrCuN10} and demonstrates the validity of our convergence theory on D-stationary points in this paper. In addition, we also show that convergence to the degenerate solution of nonlinear optimization with high accuracy is applicable for problem (TP5).

In our implementation, we use the standard starting points for $x_0$ in all test problems. The initial estimates for all Lagrange multipliers (that is, the components of $\mu_0$ and $\lambda_0$) are set to be one. {{The initial parameters are selected as}} $\rho_0=1.0$, $\sigma=0.01$, $\tau=0.5$, $\epsilon=10^{-8}$ ($10^{-6}$ in \cite{ByrCuN10}). The {{initial approximate of the Hessian $B_0$ of Lagrangian $L(x,\lambda)=\rho f(x)-\mu^Th(x)-\lambda^Tg(x)$ is simply chosen to be the identity matrix, and}} the step-size is computed by the Armijo line-search procedure. According to Definitions \ref{def0} and \ref{def1}, in order to ascertain what are the terminating points when implementing our algorithm, we calculate the measures \begin{eqnarray}
\dd\dd E_k^{\hbox{dual}}=\|\rho_{k-1}\na f(x_k)-\na h(x_k)\mu_k-\na g(x_k)\lambda_k\|_{\infty}, \nonumber\\
\dd\dd E_k^{\hbox{compl}}=\max\{\|\mu_k\circ h(x_k)+|h(x_k)|\|_{\infty}, \|\lambda_k\circ g(x_k)+\max\{0,-g(x_k)\}\|_{\infty}\}, \nonumber\\
\dd\dd E_k^{\hbox{feas}}=\max\{\|h(x_k)\|_{\infty},\|\max\{0,-g(x_k)\}\|_{\infty}\}, \nonumber\end{eqnarray}
where $\circ$ is the Hadamard product of two vectors. In addition, we will report the values of objective functions and the number of solved QP subproblems for inner iterations, the number of evaluations of functions and gradients in inner iterations.

The first test problem is referred as {\sl unique} in \cite{ByrCuN10}: \begin{eqnarray}
({\rm TP1}) \quad\min~ x_1+x_2 \quad\st~ x_2-x_1^2-1\ge 0, ~ 0.3(1-e^{x_2})\ge 0. \nonumber\end{eqnarray}
The standard initial point is $x_0=(3,2)$, an approximate DZ-stationary point close to $x^*=(0,1)$, which is also an approximate strict minimizer of the $\ell_1$-norm constraint violations, was found in \cite{ByrCuN10}. Our algorithm terminates at $k=6$ with an approximate solution to $x^*$. The output of our algorithm is given in Table \ref{tab1}, where ``iter-sb" represents the number of inner iterations for solving the $k$-th subproblem with fixed parameter $\rho_k$. \refal{alg1} takes the full step at almost all iterates. It is easy to observe from Table \ref{tab1} that both $E_{k}^{\hbox{dual}}$ and $E_k^{\hbox{compl}}$ are small enough, which show the rapid convergence to the stationary point.
\begin{table}
{\tiny
\begin{center}
\caption{Output for test problem (TP1): {$13$ inner iterations} needed.}\label{tab1} \vskip 0.2cm
\begin{tabular}{|c|c|c|c|c|c|c|c|c|}
\hline
$k$ & $f_k$ & $E_{k}^{\hbox{dual}}$ & $E_{k}^{\hbox{compl}}$ & $E_{k}^{\hbox{feas}}$ & iter-sb & $\rho_k$ & $\hbox{numf}_k$ & $\hbox{numg}_k$   \\
\hline
\hline
0 & 5 & 7 & 0 & 8 & - & 1.0 & 1 & 1  \\
1 & -0.8032 & 1.0816 & 5.0886e-09 & 1.8456 & 4 & 0.01 & 4 & 4  \\
2 & 0.9941 & 2.1018e-04 & 3.0986e-08 & 0.5155 & 5 & 1.0000e-04 & 6 & 5  \\
3 & 0.9999 & 6.1454e-08 & 2.8401e-11 & 0.5155 & 2 & 1.0000e-06 & 2 & 2  \\
4 & 1.0000 & 1.9815e-05 & 1.0122e-09 & 0.5155 & 1 & 1.0000e-09 & 1 & 1  \\
5 & 1.0000 & 3.3072e-09 & 5.6710e-11 & 0.5155 & 1 & 1.0000e-09 & 1 & 1  \\

\hline
\end{tabular}
\end{center}}
\end{table}

The second problem is the {\sl{isolated}} problem of \cite{ByrCuN10}:
\begin{eqnarray}
\min\dd\dd x\sb{1}+x_2 \nonumber\\
({\rm TP2}) \quad\quad
\st\dd\dd -x_{1}^{2}+x\sb{2}-1\ge 0, \nonumber\\
   \dd\dd -x_{1}^{2}-x\sb{2}-1\ge 0, \nonumber\\
   \dd\dd x_{1}-x\sb{2}^2-1\ge 0, \nonumber\\
   \dd\dd -x_{1}-x\sb{2}^2-1\ge 0. \nonumber
\end{eqnarray}
Starting from the same initial point {{$x_0=(3,2)$}} as problem (TP1), the algorithm in \cite{ByrCuN10} found an approximate DZ-stationary point close to $x^*=(0,0)$, a strict minimizer of the infeasibility measure in $\ell_1$ and $\ell_2$ norms. Our algorithm terminates at an approximate point to it. The output of our algorithm for problem (TP2) is reported in Table \ref{tab2}, which shows the rapid convergence to the stationary point.
\begin{table}
{\tiny
\begin{center}
\caption{Output for test problem (TP2): {$15$ inner iterations} needed.}\label{tab2} \vskip 0.2cm
\begin{tabular}{|c|c|c|c|c|c|c|c|c|}
\hline
$k$ & $f_k$ & $E_{k}^{\hbox{dual}}$ & $E_{k}^{\hbox{compl}}$ & $E_{k}^{\hbox{feas}}$ & iter-sb & $\rho_k$ & $\hbox{numf}_k$ & $\hbox{numg}_k$   \\
\hline
\hline
0 & 5 & 13 & 0 & 12 & - & 1.0 & 1 & 1  \\
1 & -0.4947 & 0.0141 & 5.0286e-12 & 1.3089 & 8 & 0.01 & 11 & 8  \\
2 & -0.0050 &  2.9454e-08 & 2.6117e-10 & 1.0025 & 4 & 1.0000e-04 & 4 & 4 \\
3 & 1.1633e-05 & 1.4189e-04 & 2.3552e-12 & 1.0000 & 1 & 1.0000e-06 & 1 & 1 \\
4 & 2.5621e-11 & 1.0276e-06 & 2.3759e-14 & 1.0000 & 1 & 1.0000e-09 & 1 & 1 \\
5 & 2.5621e-11 & 2.8630e-08 & 0 & 1.0000 & 1 & 1.0000e-09 & 0 & 0 \\

\hline
\end{tabular}
\end{center}}
\end{table}

The third test problem is the {\sl nactive} problem in \cite{ByrCuN10}:
\begin{eqnarray}
\hbox{min}\dd\dd x\sb{1} \nonumber\\
({\rm TP3}) \quad\quad
\hbox{s.t.} \dd\dd \frac{1}{2}(-x\sb{1}-x\sb{2}^2-1)\ge 0, \nonumber\\
   \dd\dd x\sb{1}-x_2^2\ge 0, \nonumber\\
   \dd\dd -x\sb{1}+x_2^2\ge 0. \nonumber
\end{eqnarray}
This problem is still infeasible. The given initial point is $x_0=(-20,10)$. The point $x^*=(0,0)$ was an infeasible stationary point with $\|\max\{0,-g(x^*)\}\|=0.5$ and is also a DL-stationary point. \refal{alg1} terminates at an approximate point $x_k$ with $k=2$ and $\rho_k=0.01$, and  takes full steps at all inner iterates. The output is shown in Table \ref{tab3} where we need to solve $12$ QP subproblems.
\begin{table}
{\tiny
\begin{center}
\caption{Output for test problem (TP3): {$12$ inner iterations} needed.}\label{tab3} \vskip 0.2cm
\begin{tabular}{|c|c|c|c|c|c|c|c|c|}
\hline
$k$ & $f_k$ & $E_{k}^{\hbox{dual}}$ & $E_{k}^{\hbox{compl}}$ & $E_{k}^{\hbox{feas}}$ & iter-sb & $\rho_k$ & $\hbox{numf}_k$ & $\hbox{numg}_k$   \\
\hline
\hline
0 & -20 & 10 & 120 & 120 & - & 1.0 & 1 & 1  \\
1 & -21.7061 & 7.9913e-11 & 1.8034e-11 & 21.7061 & 5 & 0.01 & 5 & 5 \\
2 & -4.2698e-13 & 1.4457e-10 & 1.0461e-13 & 0.5000 & 7 & 0.01 & 7 & 7 \\

\hline
\end{tabular}
\end{center}}
\end{table}

We also use \refal{alg1} to solve the challenging example given in \reff{s-example}:
\begin{eqnarray}
({\rm TP4})\quad\hbox{min}~x \quad\hbox{s.t.}~x\sp{2}-1\ge 0,~x-2\ge 0. \nonumber\end{eqnarray}
This problem has a unique global minimizer $x^*=2$, at which both LICQ and MFCQ hold, and the second-order sufficient optimality conditions are satisfied. However, like most of the methods for nonlinear optimization, including various SQP and interior-point methods, our method terminates at $x^*=-1.0000$, a DL-stationary point when starting from an infeasible point $x_0=-4$. Our numerical tests show that the point $x^*$ cannot be improved by changing the initial penalty parameter. Some new techniques should be incorporated for this difficulty.
\begin{table}
{\tiny
\begin{center}
\caption{Output for test problem (TP4): {$7$ inner iterations} needed.}\label{tab4} \vskip 0.2cm
\begin{tabular}{|c|c|c|c|c|c|c|c|c|}
\hline
$k$ & $f_k$ & $E_{k}^{\hbox{dual}}$ & $E_{k}^{\hbox{compl}}$ & $E_{k}^{\hbox{feas}}$ & iter-sb & $\rho_k$ & $\hbox{numf}_k$ & $\hbox{numg}_k$   \\
\hline
\hline
0 & -4 & 8 & 15 & 6 & - & 1.0 & 1 & 1 \\
1 & -4 & 4.1933e-10 & 7.8351e-10 & 6 & 1 & 0.1000 & 0 & 0 \\
2 & -1.0000 & 8.8707e-14 & 5.2231e-10 & 3.0000 & 6 & 0.1000 & 6 & 6 \\

\hline
\end{tabular}
\end{center}}
\end{table}

In the final experiment, we solve a standard test problem taken from \cite[Problem 13]{HocSch81}:
\begin{eqnarray}
({\rm TP5}) \quad\hbox{min}~ (x\sb{1}-2)\sp{2}+x\sb{2}\sp{2} \quad\hbox{s.t.}~ (1-x\sb{1})\sp{3}-x\sb{2}\geq 0, ~ x\sb{1}\geq 0, ~ x\sb{2}\geq 0. \nonumber\end{eqnarray}
This problem is obviously feasible, and has the optimal solution $x^*=(1,0)$ which is not a KKT point but a singular stationary point, at which the gradients of active constraints are linearly dependent. It is a challenging problem since the convergence of many algorithms depends on the LICQ. There is not much detail on the solution of this problem in the literature.

The standard initial point $x_0=(-2,-2)$ is an infeasible point.
For this problem, we select $\rho_0=10^3$. \refal{alg1} terminates at $x_{k}\approx(1.0000004967,-0.0000000000)$ with $k=26$, an approximate solution close to the minimizer $x^*$ very sufficiently. The parameter $\rho_k\to 0$ very quickly. The output of the algorithm is given in Table \ref{tab5}.
\begin{table}
{\tiny
\begin{center}
\caption{Output for test problem (TP5): {$51$ inner iterations} needed.}\label{tab5} \vskip 0.2cm
\begin{tabular}{|c|c|c|c|c|c|c|c|c|}
\hline
$k$ & $f_k$ & $E_{k}^{\hbox{dual}}$ & $E_{k}^{\hbox{compl}}$ & $E_{k}^{\hbox{feas}}$ & iter-sb & $\rho_k$ & $\hbox{numf}_k$ & $\hbox{numg}_k$   \\
\hline
\hline
0 & 20 & 18 & 29 & 2 & - & 1.0e+3 & 1 & 1 \\
1 & 20 & 7973 & 29 & 2 & 1 & 100 & 0 & 0 \\
2 & 10.0000 & 588.0000 & 9.0000 & 1.0000 & 1 & 10 & 1 & 1 \\
3 & 9.7936 & 47.6284 & 8.4966 & 1.0322 & 1 & 1 & 1 & 1 \\
4 & 3.3733 & 2.1970 & 1.1540 & 0.9907 & 3 & 0.1000 & 3 & 3 \\
5 & 0.6151 & 0.0172 & 0.0197 & 0.0098 & 11 & 1.0000e-03 & 14 & 11 \\
6 & 0.9497 & 1.0767e-07 & 1.6575e-05 & 8.2905e-06 & 8 & 1.0000e-05 & 8 & 8 \\
7 & 0.9930 & 2.0267e-06 & 3.7378e-08 & 2.3394e-08 & 6 & 3.1623e-08 & 6 & 6 \\
8 & 0.9970 & 3.2532e-06 & 4.6188e-09 & 3.0816e-09 & 2 & 5.6234e-12 & 2 & 2 \\
9 & 0.9982 & 1.1838e-06 & 1.5704e-09 & 1.0471e-09 & 1 & 1.3335e-17 & 1 & 1 \\
10 & 0.9989 & 4.6321e-07 & 3.2247e-10 & 2.1499e-10 & 1 & 4.8697e-26 & 1 & 1 \\
11 & 0.9993 & 1.7528e-07 & 7.9826e-11 & 5.3217e-11 & 1 & 1.0000e-30 & 1 & 1 \\
12 & 0.9996 & 6.7188e-08 & 1.8493e-11 & 1.2329e-11 & 1 & 1.0000e-30 & 1 & 1 \\
13 & 0.9997 & 2.5629e-08 & 4.3963e-12 & 2.9308e-12 & 1 & 1.0000e-30 & 1 & 1 \\
14 & 0.9998 & 9.7941e-09 & 1.0350e-12 & 6.9001e-13 & 1 & 1.0000e-30 & 1 & 1 \\
15 & 0.9999 & 3.7404e-09 & 2.4458e-13 & 1.6305e-13 & 1 & 1.0000e-30 & 1 & 1 \\
16 & 0.9999 & 1.4289e-09 & 1.2431e-13 & 3.8479e-14 & 1 & 1.0000e-30 & 1 & 1 \\
17 & 1.0000 & 5.4590e-10 & 8.7394e-14 & 9.0863e-15 & 1 & 1.0000e-30 & 1 & 1 \\
18 & 1.0000 & 2.0863e-10 & 6.1632e-14 & 2.1455e-15 & 1 & 1.0000e-30 & 1 & 1 \\
19 & 1.0000 & 7.9797e-11 & 4.5121e-14 & 5.0686e-16 & 1 & 1.0000e-30 & 1 & 1 \\
20 & 1.0000 & 3.0582e-11 & 3.5889e-14 & 1.1985e-16 & 1 & 1.0000e-30 & 1 & 1 \\
21 & 1.0000 & 1.1778e-11 & 3.2477e-14 & 2.8405e-17 & 1 & 1.0000e-30 & 1 & 1 \\
22 & 1.0000 & 4.5916e-12 & 3.4053e-14 & 6.7710e-18 & 1 & 1.0000e-30 & 1 & 1 \\
23 & 1.0000 & 1.8443e-12 & 4.0137e-14 & 1.6396e-18 & 1 & 1.0000e-30 & 1 & 1 \\
24 & 1.0000 & 7.9742e-13 & 4.9823e-14 & 4.1665e-19 & 1 & 1.0000e-30 & 1 & 1 \\
25 & 1.0000 & 4.1254e-13 & 6.0513e-14 & 1.2596e-19 & 1 & 1.0000e-30 & 1 & 1 \\
26 & 1.0000 & 3.0213e-13 & 6.7940e-14 & 6.4726e-20 & 1 & 1.0000e-30 & 1 & 1 \\

\hline
\end{tabular}
\end{center}}
\end{table}

In summary, the preceding numerical results not only demonstrate global convergence to the stationary points with least constraint violations on \refal{alg1} for infeasible and degenerate nonlinear programs, but also illustrate that our algorithm is capable of rapidly detecting infeasibility of nonlinear programs without driving the penalty parameter $\rho_k$ to zero for some infeasible optimization.

\section{Conclusion}

We introduce the stationary points of nonlinear optimization with least constraint violations including the D-stationary point, the DL-stationary point and the DZ-stationary point for dealing with possible infeasible optimization problems. These stationary points can be thought of as a generalization of the classic Fritz-John point, KKT point and singular stationary point for feasible optimization to infeasible optimization, where the DL-stationary point has a dependence on the objective like the KKT point, and the DZ-stationary point corresponds to the singular stationary point.
To examine the usefulness of our new stationary points, based on robustness of the exact penalty optimization, we present an exact penalty SQP method with inner and outer iterations for nonlinear optimization, and analyze its global and local convergence.
It is shown that when the solution of an infeasible optimization is a DL-stationary point, the rapid infeasibility detection can happen without driving the penalty parameter to zero. This demonstrates the commentary of Byrd, Curtis and Nocedal \cite{ByrCuN10} and is a supplement of the theory of nonlinear optimization on rapid detection of infeasibility. Some preliminary numerical results are reported.

As illustrated by the simple one-dimensional example \reff{s-example}, even for a feasible optimization problem, many state-of-the-art solvers may return an infeasible D-stationary point. It is worthwhile to design more effective numerical algorithms which can converge to a good D-stationary point.

\


\end{document}